\documentclass[a4paper,oneside,12pt]{amsart}
\usepackage[left=3cm, right=3cm]{geometry}
\usepackage{header}
\usepackage{dirtytalk}

\title[Stokes' Theorem on integral currents]{A generalized Stokes' Theorem on integral currents}
\author{Antoine Julia}
\date{\today}
\address{Département de Mathématiques, Bâtiment 307,
Faculté des Sciences d'Orsay,
Université Paris-Saclay}
\email{antoine.julia@universite-paris-saclay.fr} 
\thanks{The author is supported by the University of Padova STARS Project \say{Sub-Riemannian Geometry and Geometric Measure Theory Issues: Old and New} (SUGGESTION), and by GNAMPA of INdAM (Italy) through the project \say{Rectifiability in Carnot Groups}.}

\begin{document}
\maketitle
\begin{abstract}
  The purpose of this paper is to study the validity of Stokes' Theorem for singular submanifolds and differential forms with singularities in Euclidean space. The results are presented in the context of Lebesgue Integration, but their proofs involve techniques from gauge integration in the spirit of R.~Henstock, J.~Kurzweil and W.~F.~Pfeffer. We manage to prove a generalized Stokes' Theorem on integral currents of dimension $m$ whose singular sets have finite $m-1$ dimensional intrinsic Minkowski content. This condition applies in particular to codimension $1$ mass minimizing integral currents with smooth boundary and to semi-algebraic chains. Conversely, we give an example of integral current of dimension $2$ in $\R^3$, with only one singular point, to which our version of Stokes' Theorem does not apply.
\end{abstract}

\setcounter{tocdepth}{1}
\tableofcontents

\renewcommand\labelitemi{$\cdot$}

\section{Introduction and main results}
Stokes' Theorem is a key result in geometry and analysis. From the analytic point of view, it can be seen as a general version of the Fundamental Theorem of Calculus and of the Divergence Theorem. It lies at the core of integration by parts and thereby of the notion of weak solution of many partial differential equations. On the geometric side, Stokes' theorem is key to the De Rham Cohomology. 

 In geometric analysis, the notion of calibration (see for instance \cite{Morgan1990Calibrations, Helein1994calibration}) connects these two points of view. As for many variational problems one expects minimizers to eventually exhibit singularities, a corresponding calibration is then expected to be singular as well. It is natural to try to generalize Stokes' Theorem to singular forms and varieties. Another reason to do this is the study of PDEs on singular surfaces, in particular in order to study their weak formulations.

The classical Stokes' Theorem is usually stated as follows:
\renewcommand*{\theTheorem}{\hspace{-.15cm}}
\begin{Theorem}[Stokes' Theorem]
  If $M$ is a compact $m$-dimensional oriented $C^1$ submanifold of $\Rn$, with boundary $\partial M$ and $\omega$ is a $C^1$ differential form of degree $m-1$ on $M$, then there holds
  \begin{equation}
    \label{eq:stokes0}
    \int_M \dd \omega = \int_{\partial M} \omega.
  \end{equation}
\end{Theorem}
The main question behind the present work is the following:
\renewcommand*{\theQuestion}{\hspace{-.12cm}}
\begin{Question}
  Suppose now that $M$ is a singular submanifold with singular set $E_M$ and that $\omega$ has singularities in the set $E_\omega$. Under what conditions on the sets $E_M$ and $E_\omega$, and possibly on the types of singularities, does identity \eqref{eq:stokes0} still hold?
\end{Question}
Let us first review some classical answers to these questions. If $m=n=1$, $M$ is a compact interval and $\omega=f$ is a function, this amounts to proving a generalized Fundamental Theorem of Calculus. It is known that if $f$ is continuous, differentiable except in a countable set and $f'$ is Lebesgue integrable, then the primitive of $f'$ is equal to $f$ up to a constant. However, the Lebesgue integrability condition is not automatically satisfied, even if $f$ is differentiable everywhere --- consider for example the continous extension to $[\,-\pi^{-1/2},\pi^{-1/2}\,]$ of the function $x\mapsto x^2\sin (x^{-2})$.

 This problem was solved by introducing a new type of integral on intervals, whose first formulations were given by A. Denjoy and O. Perron. It is now known as the Henstock-Kurzweil integral \cite{Kurzweil1957generalized, Henstock1961definitions}, and can be constructed in a way very similar to the Riemann Integral, though it is more general than the Lebesgue Integral. The books \cite{Saks1964theory, Gordon1994integrals, Moonens2017integration} contain detailed presentations of these questions. The main advantage of the formulation using Riemann Sums lies in the focus on the domain of integration --- as opposed to the focus on the range, as in the Lebesgue Integral. Indeed focusing on the domain provides a better control on the behaviour of the function near pointwise singularities, by the mean of a \emph{gauge}, i.e. a non-negative function defined on the domain and controlling the size of the elements in the Riemann sums. In particular, this type of method can also yield results for the Lebesgue integral.

 If $m=n\geq 1$ and $M$ represents a bounded set of finite perimeter, a method of gauge integration was developped by W.~F. Pfeffer \cite{Pfeffer1991gauss}, following in particular works of  J. Ma\v{r}\`ik \cite{Marik1965extensions}, J. Mawhin \cite{Mawhin1981generalized}. The main result of Pfeffer Integration is a generalized Divergence Theorem, which we rephrase below as Theorem \ref{thm:Pfefferstokes} using the notation of this paper.  W.~F. Pfeffer's result extends the celebrated Theorem of E. De Giorgi and H. Federer (see \cite[Theorem~4.5.6]{FedererGMT}), which states that if $A\subset \Rm$ is a bounded set of finite perimeter and $\bv$ is a Lipschitz vector field, then the Divergence Theorem holds:
 \begin{equation*}
   \int_A \diver \bv = \int_{\partial_* A} \bv \cdot \nu_A \dd \scrH^{m-1},
 \end{equation*}
where $\partial_* A$ is the reduced boundary of $A$ and $\nu_A$ is its outer normal and for $s\in [\,0, n\,]$, $\scrH^s$ denotes the Hausdorff measure of dimension $s$ .

We now turn to the case where $m\leq n$. If $M$ is a $C^1$ submanifold, this reduces to the flat case ($m=n$) by $C^1$ triangulation and changes of variables. In order to study singular submanifolds, we choose to work in the setting of integral currents in Euclidean space. These currents were introduced in \cite{FedFlem1960} by H. Federer and W.~H. Fleming, who presented them as a satisfactory class of \say{$k$ dimensional domain of integration in euclidean $n$-space}. We will mostly follow the notation from the classical book \cite{FedererGMT}, let us introduce some of it now (see also Section \ref{sec:def}).

 Integral currents form a subset of the currents in the sense of De Rham: A current $T\in \calD_m(\Rn)$ is a continuous linear operator on $\calD^m(\Rn)$, the space of smooth differential forms of degree $m$ in $\Rn$ with compact support. The \emph{boundary} of $T$ is the current $\partial T\in \calD_{m-1}(\Rn)$ defined for $\omega \in \calD^{m-1}(\Rn)$ by $\partial T(\omega) = T(\dd \omega)$. The \emph{mass} of a current $T$ of dimension $m$ is defined as 
\begin{equation*}
    \mass(T) = \sup \{ T(\omega),\  \omega \in \calD^m(\Rn), \ \forall x\in \Rn,\ \vert \omega (x) \vert <1\}.
\end{equation*}
In particular, if $T$ represents the oriented submanifold $M$, i.e. if $T(\omega) = \int_M \omega$ for $\omega \in \calD^m(\Rn)$, then there holds $\mass(T) = \scrH^m(M)$.

For compactness purposes, one prefers to work with Lipschitz instead of $C^1$ maps. This leads to the notion of rectifiability: a set is \emph{$m$-rectifiable}, if it can be covered up to an $\scrH^m$ null set by a countable union of Lipschitz images of $\Rm$. We say that $T\in \calD_m(\Rn)$ is a \emph{rectifiable current} if it has compact support (denoted by $\spt T$) and can be represented by a triple $(M,\theta, \vect T)$, where
\begin{enumerate}[label=(\alph*)]
  \item $M\subset \Rn$ is $m$-rectifiable,
  \item  $\theta:M\to \Z$ is integrable with respect to $\scrH^m$,
  \item  $\vect T$ is an $\scrH^m$-measurable unit $m$-vector field,  $\scrH^m$ a.-e. tangent to $M$.
\end{enumerate}
We then write $T=  \Vert T\Vert \wedge \vect T$, where the Radon measure $\Vert T\Vert = \theta\scrH^m \hel M$ is called the \emph{carrying measure} of $T$; there holds $\mass(T)= \Vert T\Vert(\Rn)$.  The action of $T$ on $\omega\in \calD^m(\Rn)$ is then
\[
    T(\omega) = \int_M \langle \omega(x),\vect T(x) \rangle \dd \Vert T\Vert(x),
\]
where $\langle \cdot,\cdot \rangle : \Lambda^m(\Rn)\times \Lambda_m(\Rn) \to \R$ represents the duality pairing between $m$-covectors and $m$-vectors in $\Rn$. Finally, a current $T\in \calD_m(\Rn)$ is an \emph{integral current} ($T\in \Imrn$) if both $T$ and $\partial T$ are rectifiable currents. 

Given $T\in \Imrn$ and a smooth form $\omega\in \calD^m(\Rn)$, the identity $T(\dd \omega) = \partial T(\omega)$ can be written as
\begin{equation}
  \label{eq:stokes}
  \int \langle \dd \omega(x),\vect T (x) \rangle \dd \Vert T\Vert (x) = \int \langle \omega(x), \vect{\partial T} (x)\rangle \dd \Vert \partial T\Vert(x).
\end{equation}
In this sense Stokes' Theorem \textit{always holds} for smooth differential forms on an integral current. We are interested in the validity of \eqref{eq:stokes}, when $\omega:~\spt T\to \Lambda^m(\Rn)$ is not smooth. 
Using the theory of flat cochains,  one can also handle the case in which $\omega$ is only Lipschitz continuous  (see \cite{Whitney1957book, Federer1974flat}).

Instead of considering a global regularity condition, we can study a pointwise one: we say that a map $f: \Rm\to \R^k$ is pointwise Lipschitz at $x\in \Rn$ if there holds
\[
    \Lip_x f \defeq \lim_{r\to 0^+} \sup_{y\in \Ball(x,r)} \dfrac{\vert f(y)-f(x)\vert }{\vert y-x\vert } <+\infty.
\]
The Rademacher-Stepanov Theorem (see \cite[3.1.8 and 3.1.9]{FedererGMT}), states that if $f$ is pointwise Lipschitz at every point of a set $E\subset \Rm$, then it is differentiable Lebesgue almost everywhere in $E$.
We can now state the aforementioned theorem of W.~F. Pfeffer, in the setting of integral currents of dimension $m$ in $\Rm$.
 The statement involves the \emph{density set}, $\set_m\Vert T\Vert$ of a current $T\in \Imrn$, i.e. the subset  of $\spt T$ consisting of all the points of positive upper $m$ dimensional density of $\Vert T\Vert$ (see the definition in Section \ref{sec:def}).
\renewcommand*{\theTheorem}{\arabic{Theorem}}
\setcounter{Theorem}{-1}
\begin{Theorem}[Pfeffer, {\cite[Theorem 2.9]{pfeffer2005gauss}}]\label{thm:Pfefferstokes}
   Given $T\in \Int_m(\Rm)$ and $\omega:\spt T\to \Lambda^{m-1}(\Rm)$, identity \eqref{eq:stokes} holds provided that the following conditions are satisfied:
\begin{enumerate}[label=(\roman*)]
  \item $\omega$ is bounded in $\spt T$ and continuous in $\set_m\Vert T\Vert\backslash E_0$, where $\scrH^{m-1}(E_0)=0$,\label{thm:pfeffer:continuity}
  \item $\omega$ is pointwise Lipschitz in $\set_m\Vert T\Vert\backslash E_\sigma$, where $E_\sigma$ is $\scrH^{m-1}$ $\sigma$-finite,\label{thm:pfeffer:deriv}
  \item $\dd \omega$ is Lebesgue integrable with respect to $\Vert T\Vert$.\label{thm:pfeffer:integ}
\end{enumerate}
\end{Theorem}
 Already in dimension $1$, say on an interval, one sees that assumption \ref{thm:pfeffer:continuity} is optimal. In a way, it is also the case of condition \ref{thm:pfeffer:deriv}; indeed the nondifferentiability set of a continuous function on $[\,0,1\,]$ is of first category (see \cite{Zahorski1946ensemble}), and an uncountable set of first category contains a Cantor set (see \cite[Lemma~5.1]{Oxtoby1980MnC}), based on which one can construct a \textit{Devil's staircase} function, which does not satisfy the Fundamental Theorem of Calculus whenever the Cantor set has zero Lebesgue measure.  As for assumption \ref{thm:pfeffer:integ}, one could omit it, but at the cost of introducing another integral on bounded sets of finite perimeter, known as the Pfeffer integral (or $\calR$ integral). Following the example of \cite{pfeffer2005gauss}, we will stick to the context of Lebesgue integration while using methods coming from gauge integration.

By a smooth change of variables, Theorem \ref{thm:Pfefferstokes} could be extended to smooth oriented submanifolds with boundary. However, we are interested in singular submanifolds, so we prefer a more direct approach based on the decomposition of currents into Riemann sums. The elements of the decomposition are \emph{subcurrents}, with a definition similar but not identical to that given by E. Stepanov and E. Paolini in \cite{PaoSte2012acyclic}.

As opposed to bounded sets of finite perimeter, it is not in general possible to decompose an integral current into appropriate Riemann sums; we therefore consider currents which are locally pieces of $C^1$ submanifolds with \textit{finite perimeter}, except in a small singular set. More precisely, we say
that a current $T$ admits a \emph{$C^1$ chart} at a point $x$ of its support, if there exists a neighbourhood $U$ of $x$ such that $T\hel U$ is of the form $\phi_\# \bE_m \hel A$, where $A\subset \Rm$ is bounded and has finite perimeter, $\phi: \cl A\to \Rn$ is $C^1$ and bi-Lipschitz, and $\bE_m$ stands for the $m$-current in $\Rm$ which represents the Lebesgue measure with canonical orientation. A point in $x\in\spt T$ is called \emph{singular} $T$ does not admit a $C^1$ chart at $x$. The singular set of a current $T$ will be denoted by $E_T$ in the sequel, by definition it is a compact subset of $\spt T$.

We now introduce a way to ensure that the singular set is not \textit{too large} in a current:
A set $E\subset \spt T$ is \emph{disposable} in $T$ if there exists $C >0$ such that given $\epsilon>0$, one can find a neighbourhood $U$ of $E$ such that
\begin{equation*}
  \begin{cases}\mass(T\hel U) &<\epsilon,\\
    \mass(\partial (T\hel U)) &< C.
    \end{cases}
\end{equation*}
We say that a set $E$ is \emph{strongly disposable} if given $\epsilon$, one can find a neighbourhood $U$ of $E$ as above, such that in addition $\mass(\partial (T\hel U)) <\epsilon$. A current $T\in \Imrn$ is called \emph{weakly regular} if its singular set $E_T$ is disposable in $T$. With these notions at hand, we can now state our main result:
\renewcommand*{\theTheorem}{\Alph{Theorem}}
\renewcommand*{\theProposition}{\Alph{Theorem}}
\begin{restatable}[Generalized Stokes Theorem]{Theorem}{stokes}\label{thm:stokes}
  Let $T\in \Imrn$ be weakly regular with singular set $E_T$, and $\omega: \spt T \to \Lambda^{m-1}(\Rn)$ be a differential form satisfying
  \begin{enumerate}[label=(\roman*)]
    \item $\omega$ is bounded on $\spt T$ and continuous in $\spt T\backslash E_0$, where $E_0$ is strongly disposable in $T$,\label{thm:stokes:continuity}
    \item $\omega$ is pointwise Lipschitz on $\spt T\backslash (E_\omega\cup E_T)$, where $E_\omega$ is $\scrH^{m-1}$ $\sigma$-finite,
    \item $x\mapsto \langle \dd \omega(x), \vect T(x) \rangle$ is Lebesgue integrable with respect to $\Vert T\Vert$.
  \end{enumerate}
 Then identity \eqref{eq:stokes} holds. 
\end{restatable}
This theorem is proved in Section \ref{sec:stokes}.
An important question is whether $C^1$ charts can be replaced by bi-Lipschitz charts in the definition of weakly regular current. The only difficulty lies in establishing the validity of the derivation lemma (Lemma \ref{lemma:derivation}), as the tangent $m$-vector in a bi-Lipschitz chart is in general only $\Vert T\Vert$-approximately continuous $\Vert T\Vert$ almost-everywhere (see Remark \ref{rem:appcont}).

It is natural to ask whether the condition that $E_T$ be disposable is related to its $\scrH^{m-1}$ measure, however we show that a singular set of dimension $m-2$ can be non-disposable:
\begin{restatable}{Theorem}{contrex}\label{thm:contrex}
  There exists a current $T\in\Int_2(\R^3)$ having only one singular point and such that Stokes' Theorem in the generality of Theorem \ref{thm:stokes} does not hold for $T$, not even for continuous forms. In particular, $T$ is not weakly regular and its singular set is not disposable in $T$.
\end{restatable}
The example which we construct to prove this result can be found in Section \ref{sec:example}. 
I believe that a similar construction can yield a weakly regular current in $\Int_2(\R^3)$ with only one singular point, along with a differential form $\omega$ which is bounded on $\spt T$ and smooth on $\spt T \backslash E_T$, such that \eqref{eq:stokes} does not hold. This would show that the strong disposability condition of Theorem \ref{thm:stokes} cannot be relaxed to a Hausdorff measure condition.

The next goal is then to find practical criteria ensuring that a set is disposable. To do this, in Section \ref{sec:Minkowski}, we study the Minkowski content of a set with respect to (the carrying measure of) a current. Given a  current $T\in \Imrn$ and a set $E\subset \Rn$, the $m-1$ dimensional lower Minkowski content of $E$ with respect to $\Vert T\Vert$ is the quantity
\[
    \mathscr{M}^{m-1}_{\Vert T\Vert, *}(E) \defeq \liminf_{r\to 0} \dfrac{\Vert T\Vert (\Ball(E,r))}{r},
\]
where $\Ball(E,r)\defeq \{ x\in \Rn, \dist (x,E)< r\}$.

We first note that there are weakly regular currents whose singular set has infinite Minkowski content (Example \ref{example:infmink}); then, using classical slicing theory, we prove the following result:
\begin{restatable}{Proposition}{minkowskidispo}\label{prop:minkowskidispo}
  Given $T\in \Imrn$ and $E\subset \Rn$, the next two statements hold:
  \begin{enumerate}
    \item If $\mathscr{M}^{m-1}_{\Vert T\Vert,*} (E) < +\infty$, then $E$ is disposable in $T$, \label{minkdispo}
    \item If $\mathscr{M}^{m-1}_{\Vert T\Vert,*} (E) =0 $, then $E$ is strongly disposable in $T$.\label{minkstrongdispo}
  \end{enumerate}
\end{restatable}

This allows us to prove that codimension one mass minimizing currents without singularities at the  boundary are weakly regular (see Theorem \ref{thm:minimizing-weaklyreg}). We use a result of Cheeger and Naber \cite{CheegerNaber2013} on quantitative stratification on singular sets. These results are much stronger and also apply for example to stationary varifolds. These topics have recently attracted a lot of interest, see for instance \cite{NaberValtorta2015Varifolds}. However, note that all these results apply to conical singularities, and therefore cannot be used to study mass minimizing currents in higher codimension as these can have \textit{flat} singularities.

Other interesting classes of singular submanifolds include real algebraic and analytic subvarieties. These can be described in the framework of  o-minimal geometry, which turns out to provide enough regularity for our purpose (see the book~\cite{VDD98book} for a comprehensive presentation of o-minimal geometry). In Section \ref{sec:semi-alg} of this paper, we focus on the particular case of semi-algebraic sets (see for instance the book~
\cite{BocCosRoy1998realalg}): a semi-algebraic set in $\Rn$ is a set which can be defined by finitely many operations involving only polynomial equations and inequalities (see Definition \ref{def:semi-alg}). We define a class of integral currents called \emph{semi-algebraic chains} (see Definition \ref{def:semi-alg-chain}, these chains are also studied in \cite{KontsevichSoibelman2000Deformations, HaLaTuVo2011,Funk2016Homology}). We prove that they are weakly regular currents. Our proof is valid for chains definable in any o-minimal structure (see \cite[Chapter~6]{ThesisAJ}). A Stokes' Theorem for subanalytic varieties (which belong to the o-minimal structure of subanalytic sets) was obtained by G.~Valette in \cite{Valette2015Stokes} for \say{stratified forms}, following the works of W.~Paw\l ucki \cite{Pawlucki1985Stokes} and S.~\L ojasiewicz \cite{Lojasiewicz1991Stokes}, and a similar result for semi-algebraic varieties by L.~Shartser and G.~Valette in \cite{ShartserValette2010DeRham}. The tools we use here could lead to a more general result, allowing for singularities of the stratified forms and working in an arbitrary o-minimal structure.

To sum up, the following statement lists the classes currents for which we were able to prove that they are weakly regular.
\begin{restatable}{Theorem}{goodcurrents}\label{thm:goodcurrents}
  The following currents are weakly regular:
  \begin{enumerate}[label=(\arabic*)]
    \item Currents $T\in \Imrn$ such that $\mathscr{M}_{\Vert T\Vert,*}^{m-1} (E_T)<\infty$,\label{goodmink}
    \item Mass minimizing  integral currents of codimension $1$ with $C^{1,\alpha}$ boundary,\label{goodmini}
    \item Semi-algebraic chains.\label{goodsemialg}
  \end{enumerate}
\end{restatable}
Before moving on to the proofs, we would like to mention similar results in slightly different settings. In \cite{julia1daccepted} and \cite[Chapter~2]{ThesisAJ}, the case of one dimentional integral currents was treated, using the decomposability of these currents into sums of curves. For all $T\in \Int_1(\Rn)$, a Fundamental Theorem of Calculus is valid in the generality of Theorem \ref{thm:stokes}. Notable works on integration on domains with fractal boundaries and fractal currents include \cite{Young1934limits, HarrisonNorton1991fractalcurves, HarrisonNorton1992gauss, Harrison1999Flux, Zust2011integration, Zust2018functions}.

\subsection*{Summary of the paper}
Section \ref{sec:def} of the paper contains the preliminary definitions for most of what follows. In Section \ref{sec:HC}, we prove that a weakly regular current can be decomposed into Riemann sums in a way suitable to the proof of Stokes' Theorem, a property we call the Cousin-Howard Property.  Section \ref{sec:stokes} contains the proof of Theorem \ref{thm:stokes}, which we first reduce to the case of a continuous differential form $\omega$ (Theorem \ref{thm:stokes2}), and the construction for the proof Theorem \ref{thm:contrex}. Section \ref{sec:Minkowski} is devoted to the study of the Minkowski content condition and the proofs of Proposition \ref{prop:minkowskidispo} and of statements \ref{goodmink} and \ref{goodmini} of Theorem~\ref{thm:goodcurrents}. Finally, in Section \ref{sec:semi-alg} we present semi-algebraic chains and prove part \ref{goodsemialg} of Theorem~\ref{thm:goodcurrents}.

\subsection*{Acknowledgements}
 I wish to thank Thierry De Pauw, my PhD advisor, for introducing me to these questions and for his constant support and helpful advice.
Giovanni Alberti, Jean-Pierre Demailly and Benoît Merlet for their careful reading of my thesis and their helpful suggestions. Camillo De Lellis for his help on mass minimizing currents.  Chris Miller for an interesting exchange on the interaction between o-Minimal Geometry and Geometric Measure Theory. And Laurent Moonens for his help on the topics of Pfeffer Integration and charges. Finally, I wish to thank the anonymous referees for their thorough reading and valuable comments. I wrote my PhD thesis, on which this work is based at the Institut de Mathématiques de Jussieu, Université Paris Diderot USPC.

\setcounter{Theorem}{0}
\renewcommand*{\theTheorem}{\thesection.\arabic{Theorem}}
\renewcommand*{\theProposition}{\thesection.\arabic{Proposition}}
\section{Preliminaries}\label{sec:def}
Our notation follows mostly that of H. Federer in \cite{FedererGMT}. We work in the Euclidean space $(\Rn,\vert\cdot\vert)$. The canonical orthonormal basis is denoted $(\be_1,\dots,\be_n)$. The spaces of $m$-vectors and $m$-covectors in $\Rn$ are denoted respectively by $\Lambda_m\Rn$ and $\Lambda^m\Rn$, the norms on these spaces are also denoted by $\vert \cdot\vert$. The action of $m$-covectors on $m$ vectors is denoted by $\langle\cdot,\cdot \rangle$. The euclidean open and closed balls of center $x$ and radius $r>0$ will be denoted respectively by $\Ball(x,r)$ and $\cBall(x,r)$. Given a set $E$ in $\Rn$ and a positive $r$, the \emph{$r$-neighbourhood of $E$} is the open set $\Ball(E,r)\defeq \bigcup_{x\in E} \Ball(x,r)$. For $m=1,2,\dots$, the Lebesgue measure in $\Rm$ is denoted by $\calL^m$ and $\alpha_m$ denotes the volume of the corresponding unit ball. The $m$-dimensional Hausdorff measure in $\Rn$ is denoted by $\scrH^m$. The restriction of a measure to a set, or its multiplication by a function is denoted by $\hel$. If $\mu$ and $\nu$ are two mutually singular measures in $\Rn$, 
we write $\mu\perp\nu$. We consider $m$ dimensional integral currents in the sense of \cite{FedererGMT}: an integral current $T\in \Imrn$, of dimension $m$ in $\Rn$ can be represented by an $m$-covector valued measure, i.e.
\[
    T = \theta \scrH^m \hel M \wedge \vect T,
\]
where $M$ is a bounded $(\scrH^m,m)$-rectifiable set (see \cite[3.2.14]{FedererGMT}), $\theta$ is an integer valued $\scrH^m$-measurable function called the multiplicity of $T$, the measure $\Vert T\Vert \defeq \theta \scrH^m\hel M$ is called the carrying measure of $T$ and  $\vect T$ is a $\Vert T\Vert$-measurable field of unit length $m$-vectors tangent to $M$ at $\Vert T\Vert$-almost every point. Integral currents of top dimension: $m$ in $\Rm$ with multiplicity one and positive orientation are special: such a current is of the form $\bE^m \hel A$, where $\bE^m \defeq \calL^m \wedge \be_1\wedge \dots \wedge \be_m$ and $A$ is a bounded set of finite perimeter in $\Rm$ (see for instance \cite{AmFuPa2000} or \cite[Section 4.5]{FedererGMT}).
 The support of a current $T$ is denoted by $\spt T$ (note that $\spt T= \spt \Vert T\Vert$). A current $T$ has boundary $\partial T$, mass $\mass(T)$ and flat norm $\flatn(T)$.

Given a measure $\mu$ on $\Rn$, and $k\in \{0,\dots,n\}$, the $k$-dimensional upper-density of $\mu$ at a point $x\in \Rn$ is given by
\[
     \Theta^{k,*}(\mu,x) \defeq \limsup_{r\to 0} \dfrac{\mu(\cBall(x,r))}{\alpha_k r^k}.
\]
The set of positive $k$-dimensional upper density points of the measure $\mu$ is denoted by $\set_k \mu$.

 The \emph{essential closure} of a set $A$ in $\Rm$ is defined by 
\[
\cl_e A \defeq \{x\in \Rm, \Theta^{m*}(\calL^m\hel A,x)>0\} = \set_m (\calL^m\hel A),
\]
 it coincides with $\set_m \Vert \bE^m\hel A\Vert $. Note that $\cl_e A$ is contained in the topological closure $\cl A$ of $A$ but that $\cl A\backslash \cl_e A$ can have positive Lebesgue measure. More importantly, in positive codimension, there exists $T\in \Imrn$ such that $\spt T\backslash \set_m \Vert T\Vert$  $\spt T$ has positive $\calL^n$ measure. However there always holds $\calL^n(\cl_e A)= \calL^n(A)$ and $\mass(T) = \Vert T\Vert(\set_m\Vert T\Vert) =\Vert T\Vert(\Rn)$.

Finally, 
 given a current $T\in \Imrn$ with $1\leq m\leq n$ and a Lipschitz function $f:\Rn\to \R$, for $r\in \R$, the \emph{slice of $T$ by $f$ at $r$} is defined as
\[
    \langle T, f,r\rangle \defeq (\partial T)\hel \{x,f(x)>r\}-\partial(T\hel \{x, f(x)>r\}),
\]
and for almost all $r\in \R$, $\langle T,f,r\rangle$ is an integral current of dimension $m-1$ (see \cite[4.2.1]{FedererGMT}). Furthermore, the mass of the slices is controlled by the total mass and the Lipschitz constant of $f$ as follows:
\begin{equation}\label{eq:control-slices}
    \int_{-\infty}^{+\infty} \mass (\langle T,f,r\rangle) \dd r \leq \Lip(f) \mass(T).
  \end{equation} 
  In this paper, we will consider slicing functions $f$ of the type $\dist(\cdot, E)$, where $E$ is a subset of $\spt T$, sometimes containing a single point.

We now turn to less classical concepts, which will be used in the sequel.
\begin{Definition}\label{def:chart}
 The current $T\in \Imrn$ has a \emph{$C^1$ chart} in an open set $V\subseteq \Rn$ if there exists an integer $\theta$, a bounded set of finite perimeter $A\subset \Rm$ and a bi-Lipschitz $C^1$ map $\phi:\cl A\to V$, with 
\[
     T\hel V = \theta \phi_\# (\bE^m\hel A).
\]
 We identify such a $C^1$ chart with the $4$--uple $(\theta,V ,A,\phi)$. 
Note that this notion of chart is weaker that the one used to define differentiable manifolds, in that it allows the current to have many holes, like a set of finite perimeter.
A point $x\in \spt T$ is a \emph{regular point} if $T$ admits a $C^1$ chart as above with $x\in V$.
\end{Definition}
Our aim is to derive integration results on currents by the use of Riemann sums. We therefore need a way to decompose currents into \textit{small pieces}. Given $T\in \Imrn$ and a $\Vert T\Vert$ measurable set $A\subseteq \Rn$, recall that $T\hel A$ is a rectifiable current. If $S\defeq T\hel A$ is integral, we say that $S$ is a \emph{subcurrent} of $T$ and we write $S\sqsubset T$.  We denote the space of subcurrents of $T$ by $\subcspace (T)$. In the present paper a few simple properties of $\subcspace(T)$ will be used, which we list here. The proofs of these statements can be found in \cite[Section 3.1]{ThesisAJ},  where a more thorough study was carried out.
\begin{Proposition}
  For $T\in \Imrn$, the following statements hold:
  \begin{enumerate}[label=(\roman*)]
    \item If $\theta$ is a non-zero integer, then $\subcspace(\theta T) = \theta \subcspace(T)$.\label{subc multiple}
    \item For $S \sqsubset T$, there holds $\vect S =\vect T$ $\Vert S\Vert$ almost everywhere. Furthermore, we have $\set_m \Vert S\Vert \subseteq \set_m \Vert T\Vert$ and $\spt S\subseteq \spt T$.\label{subc density}
    \item $S \in \Imrn$ is in $\subcspace (T)$ if and only if $T-S\in \subcspace(T)$, and if and only if $\Vert T-S\Vert \perp \Vert S\Vert$.\label{subc perp}
    \item If $S\sqsubset T$ and $R\sqsubset S$, then $R\sqsubset T$. \label{subc transitivity}
    \item If $S$ and $S'$ are subcurrents of $T$ with $\Vert S\Vert \perp \Vert S'\Vert$, then $S+ S'$ is a subcurrent of $T$.\label{subc addition}
    \item If $\phi : \spt T\to \R^{n'}$ is bi-Lipschitz, then $\phi_\# \subcspace(T) = \subcspace(\phi_\# T)$. \label{subc pushforward} Lipschitz continuity of $\phi$ alone is not sufficient. 
    \item If $S\sqsubset T$ and $f: \Rn \to \R^k$ is Lipschitz, then for $\calL^k$ almost all $y\in \R^k$,
\[
     \langle S,f,y \rangle \sqsubset \langle T,f,y\rangle.
\]
\label{subc slice}
  \end{enumerate}
\end{Proposition}

We can now study functions on the space $\subcspace(T)$.  Given such a function, $F$, we say that it is \emph{additive} 
if whenever $S$ and $S'$ are non-overlapping subcurrents of $T$, $F(S+S') = F(S)+F(S')$.
(Two currents $S$ and $S'$ are called \emph{non-overlapping} if their carrying measures are mutually singular: $\Vert S \Vert \perp \Vert S'\Vert$.) We say that $F$ is \emph{continuous} if for all sequence $(S_j)_j$ in $\subcspace(T)$ with 
\begin{equation*}
  \begin{cases}
   \sup_j (\mass( S_j)+  \mass(\partial S_j))<+\infty  \\
  \text{ and }\quad  \flatn(S_j)\to 0,
 \end{cases}
\end{equation*}
 then $F(S_j) \to 0$ as $j$ tends to infinity. 

\begin{Example}\label{example:additive}
One of the additive functions we are most interested in is the \textit{circulation} (or \textit{rotation}) of a continuous $m-1$ form $\omega$ defined on $\spt T$, which we denote by $\Theta_\omega$ it is defined for $S\in\subcspace(T)$ by
\[
    \Theta_\omega(S) \defeq \int \langle \omega(x), \vect{\partial S}(x) \rangle \dd \Vert \partial S\Vert(x).
\]
In particular, if $\omega$ is a smooth differential form with compact support, then there holds $\Theta_\omega(S) =\partial S(\omega)$.
For a general continuous form $\omega$, let us prove the additivity and continuity of $\Theta_\omega$. If $S$ and $S'$ are two subcurrents, $\partial S+ \partial S' = \partial( S+ S')$, and additivity is clear. 

To see that $\Theta_\omega$ is continuous, fix $\epsilon>0$ and consider a smooth $(m-1)$ form $\omega_\epsilon$ with $\vert \omega-\omega_\epsilon\vert_\infty<\epsilon$. Given a sequence $(S_j)_j$ converging to $0$ in the flat norm with uniformly bounded mass and boundary mass, we also have $\flatn(S_j)\to 0$, thus $\Theta_{\omega_\epsilon}(S_j) \to 0$ by definition of flat convergence. Furthermore, for all $S\in \subcspace(T)$,
\begin{equation*}
  \vert \Theta_\omega(S) -\Theta_{\omega_\epsilon}(S)\vert \leq \int \langle \omega(x)-\omega_\epsilon(x), \vect{\partial S}(x) \rangle \dd \Vert \partial S\Vert(x)\leq \epsilon \mass(\partial S).
\end{equation*}
Thus for $j$ large enough, $\vert \Theta_{\omega}(S_j) \vert\leq \epsilon(1+\mass(\partial S_j)) \leq C\epsilon.$ As $\epsilon$ is arbitrary, $\Theta_\omega(S_j)\to 0$ as $j\to \infty$. Finally, notice that if $\omega$ is only a bounded Borel measurable form (say defined $\scrH^{m-1}$ almost-everywhere on $\spt T$), $\Theta_\omega$ is still additive, but not necessarily continuous.
\end{Example}
The following elementary facts will be useful:
\begin{Proposition}
For $T\in \Imrn$.
  \begin{enumerate}
    \item Continuous additive functions on $\subcspace(T)$ form a vector space. 
    \item If $\phi$ is a bi-Lipschitz map from $\spt T$ to a subset of $\R^{n'}$, the push-forward of $T$ by $\phi$: $T'\defeq \phi_\# T$ is integral and whenever $G$ is a continuous additive 
function on $\subcspace(T')$, we can define the \emph{pullback} $F=\phi^\# G$ of $G$by $\phi$ as
\[
     \forall S\in \subcspace(T), \quad F(S) = (\phi^\# G) (S) \defeq G
(\phi_\# S).
\]
 $F$ is  a continuous additive 
function on $\subcspace(T)$.
    \item Similarly if $F$ is a function on $\subcspace(T)$, and $\theta$ is a nonzero integer, then $\theta F$ is a function on $\subcspace(\theta T)$. Continuity and additivity are preserved. 
  \end{enumerate}
\end{Proposition}

The next result is essential for integration purposes:

\begin{Proposition}
  Given $T\in \Imrn$, the mass operator $\mass$ restricted to $\subcspace(T)$ is additive and continuous.
\end{Proposition}
The continuity part of the statement might be surprising, as mass is usually only lower-semi continuous for the convergence of integral currents.  Compare with the case of the sequence $S_j \defeq \lseg (0,0), (j^{-1},1)\rseg + \lseg (j^{-1},1), (2j^{-1},0)\rseg$ in for $j=1,2,\dots$, where for $a,b\in \R^2$, $\lseg a,b\rseg$ is the $1$ current associated to the oriented segment joining $a$ to $b$. The limit of the $S_j$ is $0$; yet they all have mass larger than $2$. However, the $S_j$ cannot all be subcurrents of the same integral current. Let us pass to the proof.
\begin{proof}
Additivity is clear. To prove the continuity of $\mass\vert_{\subcspace(T)}$, consider a sequence of subcurrents of $T$: $(S_j)_{j=1,2,\dots}$ such that $\sup_j \mass(\partial S_j) <+\infty$ and $S_j\to 0$ in the flat norm. Clearly $\mass(T) = \mass(T-S_j) + \mass(S_j)$ and $T-S_j$ tends to $T$ in the flat norm. By lower semi-continuity of the mass, there holds
   \begin{equation*}
     \mass(T) \leq \liminf_{j\to \infty} \mass(T-S_j) = \mass(T)-\limsup_{j\to \infty} \mass(S_j).
   \end{equation*}
Thus $\limsup_j \mass(S_j) \leq 0$ and $\mass(S_j)$ tends to $0$.
\end{proof}
Note that, as a consequence, whenever $\zeta$ is a $\Vert T\Vert$ essentially bounded $m$-form on $\spt T$, the additive function
\[
   \Xi_\zeta : S\in \subcspace(T) \mapsto \int \langle \zeta(x), \vect S (x) \rangle \dd \Vert S\Vert(x)
\]
is also continuous, as $\vert \Xi\zeta(S)\vert \leq \Vert \zeta\Vert_{\mathrm{L}^\infty} \mass(S)$ for all $S\in \subcspace(T)$.
\section{The Cousin-Howard Property}\label{sec:HC}
 Given $T\in \Imrn$, we call \emph{tagged family in $T$} any collection $\calP$ of pairs $(x,S)$ such that
\begin{enumerate}[label =(\alph*)]
  \item $\forall (x,S)\in \calP$, $S\neq 0$ is a subcurrent of $T$ and $x\in \spt S$,
  \item $\forall (x,S), (x',S') \in \calP$, either $(x,S)=(x',S')$ or $S$ and $S'$ are non-overlapping.
\end{enumerate}
If for all $(x,S)\in \calP$ the point $x$ is contained in a set $A\subset \spt T$, we say that $\calP$ is \emph{based in $A$}. The \emph{body} of the tagged family $\calP$ is the subcurrent $[\,\calP\,] \defeq \sum_{(x,S)\in\calP} S$. Given a function $\delta : A\subset \spt T\to [\,0, +\infty )$, $\calP$ is \emph{$\delta$-fine} if all $(x,S)\in \calP$ satisfy $x\in A$ and $\diam \spt S <\delta(x)$. Such a $\delta$ is called a \emph{gauge} in $A$ if $\{\delta = 0\}$ is $\scrH^{m-1}$ $\sigma$-finite.

Given a number $\eta>0$, a non-zero current $S\in \Imrn$ is called $\eta$-\emph{regular} if there holds
\[
   \text{reg}(S) \defeq \dfrac{\mass(S)}{\mass(\partial S)\diam \spt S} >\eta.
\]
If $S$ is non-zero and verifies $\partial S=0$, we will set $\text{reg} (S)= +\infty$. If $S\in \Int_m(\Rm)$ represents a cube, it is clear that $\text{reg}(S) =(2 m\sqrt{m})^{-1}$.
The number $\reg(S)$ is called the \emph{regularity of $S$}. If $\eta : A\subset \spt T\to \Rn$ is a nonnegative function, we say that the tagged family $\calP$ is \emph{$\eta$-regular} if for all $(x,S)\in \calP$, $x\in A$ and $\reg(S) >\eta(x)>0$. In particular any $\eta$-regular $\calP$ is necessarily based  in $\{x,\eta(x)>0\}$. Given a function $G$ on $\subcspace(T)$ and a number $\epsilon>0$, we say that a tagged family $\calP$ in $T$ is \emph{$(G,\epsilon)$-full (in $T$)} if there holds
\[
    \vert G(T-[\,\calP\,]) \vert <\epsilon.
\]

\begin{Definition}
  An integral currrent $T\in \Imrn$ has the \emph{Cousin-Howard Property} if there exists a function $\eta_T$ defined on $\spt T$, such that given
  \begin{enumerate}[label=(\alph*)]
    \item a function $\eta:\spt T\to \R^+$ satisfying for $x\in \spt T$
 \begin{equation*}
     \begin{cases}
               0< \eta(x) <\eta_T(x) &\text{ if } \eta_T(x) >0,\\
               \eta(x)=0 &\text{ if } \eta_T(x) = 0,
               \end{cases}
  \end{equation*}
    \item a gauge $\delta$ defined on $\spt T\backslash E_T$,
    \item a finite collection of continuous additive functions $G_1,\dots,G_p$ on $\subcspace(T)$,
    \item a positive real number $\epsilon$,
  \end{enumerate}
then there exists a $\delta$-fine, $\eta$-regular tagged family in $T$, which is $(G_j,\epsilon)$-full for $j=1,\dots,p$.
\end{Definition}
Recall that a current is weakly regular  if its singular set $E_T$ is \emph{disposable}, that is, if there exists a positive constant $C$ such that for every positive $\epsilon$, one can find a subcurrent $T_\epsilon$ of $T$, for which there holds
\begin{equation*}
  \begin{cases}
    \spt T_\epsilon \cap E_T&=\emptyset,\\
    \mass(T-T_\epsilon) &<\epsilon,\\
    \mass(\partial(T-T_\epsilon))&<C.
  \end{cases}
\end{equation*}

The aim of this section is to prove the following result:
\begin{Theorem}\label{thm:weaklyreg-hc}
  If $T$ is weakly regular then it has the Cousin-Howard Property for some function $\eta_T$, positive outside of the singular set of $T$.
\end{Theorem}
To prove this, we need two intermediate results: first, we recall the Cousin-Howard Lemma, which states that bounded sets of finite perimeter have the Cousin-Howard Property (Theorem \ref{thm:HCLemma}). We then extend this to bi-Lipschitz pushforwards of sets of finite perimeter (Lemma \ref{lemma:bilipHC}). Finally, we consider a weakly regular current $T$ and define a function $\eta_T$ on $\spt T$, depending on a choice of charts covering the regular set of $T$. Using the previous Lemma  we prove that $T$ has the Cousin-Howard Property for $\eta_T$.

\subsection{Bounded sets of finite perimeter}
Bounded sets of finite perimeter in $\Rm$ can be represented by integral currents of dimension $m$.  It is equivalent to consider a bounded set $A$ of finite perimeter in $\Rm$ and the corresponding integral current $\bE_m\hel A$. The subcurrents of $\bE_m\hel A$ correspond to the subsets of finite perimeter of $A$ and for such a set $B\subset A$, there holds
\begin{equation*}
  \begin{cases}
   \mass(\bE_m\hel B) = \calL^m(B),\\
   \mass(\partial (\bE_m\hel B)) = \mathrm{Perimeter}(B),\\
   \diam \spt (\bE_m\hel B))= \diam \cl_e B,
\end{cases}
\end{equation*}
where $\cl_e B = \set_m \Vert \bE_m\hel B\Vert$ is the essential closure of $B$.
Thus, tagged families, fineness and regularity can be defined indifferently with subsets of finite perimeter or with subcurrents.

\begin{Theorem}[Cousin-Howard Lemma]\label{thm:HCLemma}
  Let $A$ be a bounded set of finite perimeter in $\Rm$. Given a finite collection $G_1,\dots,G_p$ of continuous additive function  on $\subcspace(\bE_m\hel A)$, $\epsilon>0$, a function $\eta$ with values in $(0,2^{-1} m^{-3/2})$ and a gauge $\delta$ on $\cl_e A$, there exists a $\delta$-fine, $\eta$-regular tagged family in $A$ which is $(G_j,\epsilon)$-full for each $j\in \{1,\dots,p\}$. In particular, $\bE^m\hel A$ has the Cousin-Howard Property with the function $\eta_T= 2^{-1}m^{-3/2}$ on $\cl_e A$.
\end{Theorem}
We only sketch the proof, which can be found in \cite[section 2.6]{Pfeffer2001book} or in \cite[chapter~4]{ThesisAJ}. Although, in the former, $\eta$ is always a constant, no crucial change is required in order to allow $\eta$ to be a function.
\begin{proof}[Sketch of proof]
 Recall that for $r>0$, the cube $[\,0,r\,]^m$ is a bounded set of finite perimeter with regularity $2^{-1}m^{-3/2}$. Thus, if we can decompose our set into cubes we are done. In the case where $A$ is a $1$-dimensional interval, and $\delta$ is positive, the existence of a decomposition into a tagged family of dyadic intervals is known as Cousin's Lemma and proved by contradiction. This extends to the case of cubes in any dimension. The tagged family then covers $A$ completely in the sense that $[\calP]=\bE^m\hel A$.

 Then let $F = \max\{\vert G_j\vert, j=1,\dots,p\}$. If $\delta$ is allowed to take the value zero in an $\scrH^{m-1}$ $\sigma$-finite set, one uses the equivalence of the Hausdorff measure and the net measure (see Falconer \cite[Theorem 5.1]{Falconer1985geometry}) and the continuity of $F$ to define a positive gauge $\tilde{\delta}$ on $A$, by modifying $\delta$ on $\{\delta = 0\}$. Then using Cousin's Lemma one gets a $\tilde{\delta}$-fine tagged family which covers $A$. If $\tilde{\delta}$ is defined correctly, the cubes tagged at a point in $\{\delta=0\}$ have a small contribution to $F$ and they can be removed, leaving us with a $\delta$-fine tagged family.

In order to consider sets which are not cubes, one needs to look at what happens at the boundary. If $A$ is essentially closed --- meaning that $\cl_e A$ is closed --- one can inscribe $A$ in a cube $Q_0$ and define a gauge $\tilde{\delta}$ so that $\tilde{\delta}$-fine sets tagged at points outside of $\cl A$ do not intersect $A$. It is also possible to choose $\tilde{\delta}$ so that if $Q$ is $\tilde{\delta}(x)$-fine with $x\in \cl_e A\backslash \set_{m-1}(\scrH^{m-1}\hel \partial_* A)$ and regular, then $Q\cap \cl_e A$ is also regular. Note that  $\scrH^{m-1}(\set_{m-1}(\scrH^{m-1} \hel \partial_*))<+\infty$. We then apply Cousin's lemma to the cube $Q_0$ and the gauge $\tilde{\delta}$. From a $\tilde{\delta}$ fine tagged family composed of cubes $Q$, one gets a tagged family in $A$ whose elements are of the form $(x,Q\cap \cl_e A)$. 

Lastly, if $A$ is not essentially closed, a result of Tamanini and Giacomelli \cite{TamGia1990Monotone} ensures that $A$ can be approximated from the inside by an essentially closed subset of finite perimeter, controlling the perimeter of the difference between $A$ and the approximating sets, which is sufficient thanks to the continuity of $F$. (In the language of this paper, the result of \cite{TamGia1990Monotone} states that $\cl A\backslash \cl_e A$ is disposable for $\bE^m\hel A$.)
\end{proof}

In the sequel, we pass to higher codimension and use the above result inside $C^1$ charts of currents. Let us first show that the Cousin-Howard Property is stable under multiplication by an integer and under bi-Lipschitz push-forwards:

\begin{Lemma}\label{lemma:bilipHC}
  If an integral current $T\in \Imrn$ is of the form:
\[
    T = \theta \phi_\# (\bE^m\hel A),
\]
where $A$ is a bounded set of finite perimeter in $\Rm$, $\phi:A\to \Rn$ is bi-Lipschitz and $\theta$ is a non zero integer, then $T$ has the Cousin-Howard Property with the function $\eta_T$ verifying for all $x\in \spt T$:
\[
    \eta_T(x) \defeq \left (\dfrac{1}{\Lip \phi \Lip (\phi^{-1})}\right)^m \dfrac{1}{2m^{3/2}}.
\] 
\end{Lemma}

\begin{Remark}\label{obs:bilip}
  Before proving this lemma, let us justify the expression for $\eta_T$. Given
  \begin{enumerate}[label=(\alph*)]
    \item an integral current $R\in \Imrn$, with regularity, $\reg R>\eta>0$,
    \item a bi-Lipschitz function $\phi : \spt R\to \R^{n'}$,
    \item a non-zero integer $\theta$,
  \end{enumerate}
  then the push-forward current $\theta \phi_\# R$ is in $\Int_m(\R^{n'})$ and satisfies
  \begin{eqnarray*}
   \mass(\theta \phi_\# R) &\geq& \theta (\Lip \phi^{-1})^{-m} \mass(R),\\
    \mass(\theta \phi_\#R)  &\leq& \theta (\Lip \phi)^{m-1}\mass(\partial R), \\
    \diam (\spt(\theta \phi_\# R)) &\leq&  (\Lip \phi) \diam \spt R.
  \end{eqnarray*}
  Thus there holds
\[
     \reg(\theta \phi_\# R) = \dfrac{\theta \mass(\phi_\# R)}{\theta \mass(\partial(\phi_\#R)) \diam (\spt\phi_\# R)} > \dfrac{1}{(\Lip \phi \Lip (\phi^{-1}))^{m}}\eta.
\]
\end{Remark}
\begin{proof}[Proof of Lemma \ref{lemma:bilipHC}]
  Fix a gauge $\delta$ on $\set_m\Vert T\Vert$, a finite collection of continuous additive functions $G_1,\dots,G_p$ on $\subcspace(T)$ and a positive number $\epsilon$. We define a gauge $\tilde{\delta}$ on $A$ by
\[
    \forall x\in \set_m(\bE_m\hel A),\, \tilde{\delta}(x) \defeq \delta(\phi(x)) /\Lip(\phi)
\]
and $\tilde{\delta}(x)\defeq 0$ for $x \in A\backslash \set_m(\bE_m\hel A)$. The pullback functions $\phi^\# G_j$ defined for $j=1,\dots,p$ by
\[
   \forall S\in \subcspace(\bE_m\hel A),\, (\phi^\#G_j)(S) \defeq G_j(\phi_\# S),
\] 
are continuous and additive, so are the $\theta (\phi^\# G_j)$. Notice that $\subcspace(\bE^m\hel A)$ corresponds exactly to the family of currents in $\mathbb{I}_m(\Rm)$ that represent a bounded subset of $A$ with finite perimeter. Given $\eta\in (0,2^{-1}m^{-3/2})$ we can apply the Cousin-Howard Lemma to $(\bE^m\hel A)$, the gauge $\tilde{\delta}$, the functions $\phi^\# G_j$ and the error term $\epsilon/\vert \theta\vert$. This yields a $\tilde{\delta}$-fine $\eta$-regular tagged family $\calP$ in $\bE^m\hel A$ such that for all $j=1,\dots,p$ there holds
\[
     (\phi^\#G_j)(\bE^m\hel A -[\calP]) <\dfrac{\epsilon}{\vert \theta\vert},
\]
and members of $\calP$ are of the form $(x,\bE^m\hel B)$, where $B$ is a subset of $A$ of finite perimeter and $x\in \cl B$. 

For such a member of $\calP$ let $S\defeq \theta \phi_\# (\bE^m\hel B)$. Note that $S$ is a subcurrent of $T$ with $\phi(x)\in \spt S$; since $\phi$ is bi-Lipschitz, the collection $\calP' \defeq \{(\phi(x),\theta \phi_\# S), (x,S)\in \calP\}$ is a tagged family in $T$. Furthermore, $\calP'$ is $\delta$-fine and applying Observation \ref{obs:bilip} we infer that $\calP'$ is $\eta'$-regular, where $\eta' \defeq (\Lip \phi \Lip(\phi^{-1}))^m \eta$. Finally, for $j=1,\dots, p$, $\calP'$ also satisfies:
\[
    G_j\left (T-\sum_{(y,S)\in \calP'} S\right ) =\theta (\phi^\# G_j)\left (\bE^m\hel A - \sum_{(x,\bE^m\hel B)\in \calP} \bE^m\hel B \right ) < \epsilon.
\]
This proves that  $T=\theta \phi_\# (\bE_m\hel A)$ has the Cousin-Howard Property for $\eta_T$.
\end{proof}
\subsection{Disposable sets and the Cousin-Howard Property}
 With the language of subcurrents, we can give an equivalent definition of disposability:
A set $E\subset \Rn$ is \emph{disposable} in an integral current $T\in \Imrn$ if there exists $C>0$ such that for all $\epsilon>0$ one can find a subcurrent $T_\epsilon$ of $T$ with the following properties:
\begin{equation}\label{eq:subcdispo}
  \begin{cases}
    \spt T_\epsilon \cap E&=\emptyset,\\
    \mass(T-T_\epsilon) &<\epsilon,\\
    \mass(\partial(T-T_\epsilon))&<C.
  \end{cases}
\end{equation}
A current whose singular set is disposable is called \emph{weakly regular}. 
We can now prove the main result of this section:
\begin{proof}[Proof of Theorem  \ref{thm:weaklyreg-hc}]
  We start by choosing a locally finite cover of the regular set of $T$ by $C^1$ charts along with a suitable function $\eta_T$. This choice is not unique as it relies on a (para)compactness argument. Without loss of generality, we can choose the $C^1$ charts at regular points of $T$ to be based on open balls. As $\spt T\backslash E_T$ is paracompact we can consider a countable, locally finite cover associated to a countable collection of regular points $(x_j)_j$. Denote these charts by $(\theta_j,A_j,\phi_j,\Ball(x_j,r_j))$.

For $y\in \spt T\backslash E_T$, there are finitely many indices $j$ such that $y\in \Ball(x_j,r_j)$, corresponding to finitely many charts. We choose $\eta_T(y)$ corresponding to the least regular of these chart; more precisely let
\[
   \eta_T(y)\defeq \min\{ (\Lip \phi_j \Lip(\phi_{j}^{-1}))^{-m} m ^{-3/2}/2, y\in\Ball(x_j,r_j)\}.
\]
Let $j(y)$ be the integer corresponding to a realization of this minimum:
\begin{equation}\label{eq:choice-of-index}
   j(y) = \text{argmin} \{(\Lip \phi_{j} \Lip(\phi_{j}^{-1}))^{-m}, y\in \Ball(x_j,r_{j})\}.
\end{equation}
 We now choose the chart at $y$ to be a restriction of the chart at $x_{j(y)}$,  picking $r_y$ such that $\Ball(y,r_y)\subset \Ball(x_{j(y)},r_{j(y)})$ and $T\hel \Ball(y,r_y)$ is integral.  Let the chart of $T$ at $y$ be
\[ 
   (\theta_y,A_y,\phi_y, \Ball(y,r_y)) \defeq \left (\theta_{j(y)}, A_{j(y)}\cap \phi_{j(y)}^{-1}(\Ball(y,r_y)), \phi_{j(y) }\vert_{ \phi_{j(y)}^{-1}(\Ball(y,r_y))}, \Ball(y,r_y)\right).
\]
The reason we chose the charts and $\eta_T$ in such a way, is that now for $y\in \spt T \backslash E_T$ and $y'\in \Ball(y,r_y)$, we have $y'\in \Ball(x_{j(y)},r_{j(y)})$ and thus
\begin{equation}
  \label{eq:reg-at-x}
  \eta_T(y') \leq \eta_T(y).
\end{equation}

Let us now prove that $T$ has the Cousin-Howard Property for the function $\eta_T$.
 Pick $\epsilon>0$ and continuous additive functions $G_1,\dots G_p$ on $\subcspace (T)$ as well as a function $\eta: \spt T\to \R$ such that for $y\in\spt T\backslash E_T$: $0<\eta(y)  <\eta_T(y)$ and $\eta(E_T)=\{0\}$.
  \begin{Claim}\label{claim subc ET}
There exists a subcurrent $T_\epsilon$ of $T$ with $E_T\cap \spt T_\epsilon =\emptyset$ and $G(T-T_\epsilon) <\epsilon/2$.
  \end{Claim}
  \begin{innerproof}[Proof of Claim \ref{claim subc ET}]
    Let $C$ be the constant associated to the definition of disposability of $E_T$ in $T$. By the continuity of $G$, there exists $\tau>0$ such that whenever $S\sqsubset T$ satisfies $\mass(S) <\tau$ and $\mass(\partial S) \leq  \mass(\partial T) + C$, there holds $\vert G(S)\vert \leq \epsilon/2$. By the disposability of $E_T$, there exists a subcurrent $T_\epsilon$ of $T$ such that the conditions of \eqref{eq:subcdispo} hold.
  Thus $\mass(\partial(T-T_\epsilon) )<C+\mass(\partial T)$; therefore $G(T-T_\epsilon) <\epsilon/2$, as claimed.
  \end{innerproof}
Such a $T_\epsilon$ being fixed, note that $\spt T_\epsilon$ is compact. Consider the charts of $T$ $(\theta_y,A_y,\phi_y,\Ball(y,r_y))$ as above. The collection of open balls $(\Ball(y,r_y/2))_{y\in \spt T_\epsilon}$ covers $\spt T_\epsilon$ and we can extract a finite subcover, say associated to the points $y_1,\dots,y_q$.

As $T_\epsilon \in \Imrn$, by classical slicing theory (see e.g.~\cite[4.2.1]{FedererGMT}), we can pick $r_1\in (r_{y_1}/2,r_{y_1})$ such that $T_\epsilon \hel \Ball(x_1,r_1)\sqsubset T_\epsilon$. Let $T_1 \defeq T_\epsilon \hel \Ball(y_1,r_1)$. Note that $T_\epsilon - T_1\in \Imrn$. We can thus repeat the argument: for $j\in \{2,\dots, q-1\},$ fix $r_j\in (r_{y_j}/2, r_{y_j})$ such that
\[
    T_j\defeq \left( T_\epsilon -\sum_{k=1}^{j-1} T_j\right ) \hel \Ball(y_j,r_j)\sqsubset   \left (T_\epsilon -\sum_{k=1}^{j-1} T_j\right ).
\]
Finally, we let $T_q \defeq T_\epsilon -\sum_{j=1}^{q-1} T_j$. It can be that for some $j$ the current $T_j = 0$, we avoid this by relabeling the sequence and taking $q$ smaller. The $T_j$ form a collection of pairwise non-overlapping subcurrents of $T_\epsilon$, with $T_\epsilon = \sum_{j=1}^q T_j$. Each $T_j$ is supported inside the ball $\Ball(y_j,r_j)$ and is of the form:
\[
    T_j \defeq \theta_{y_j} \phi_{y_j \#} (\bE^m \hel A'_j),
\]
where $A'_j$ is a subset of finite perimeter of $A_{y_j}$. Notice also that by the choice of $\eta_T$, given $x\in \set_m\Vert T_j\Vert$, we have  $x\in\Ball(y_j,r_{y_j})$ and by \eqref{eq:reg-at-x} there holds
\[
   \eta(x) < \eta_T(x) \leq  (\Lip \phi_{y_j}\Lip(\phi_{y_j}^{-1}))^{-m} m ^{-3/2}/2\eqdef\eta_j.
\]

For $j=1,\dots,q$, we use the fact that $T_j$ has the Cousin-Howard Property (Lemma \ref{lemma:bilipHC}) and apply it with the subadditive function $G\vert_{\subcspace(T_j)}$, the gauge $\delta\vert_{\set_m\Vert T_j\Vert}$ and the error $\epsilon/(2q)$ to get an $\eta_j$-regular $\delta$-fine tagged family $\calP_j$ in $T_j$ such that for $i=1,\dots,p$
\[
    G_i(T_j-[\calP_j]) < \epsilon/(2q).
\]
Concatenating the tagged families $\calP_j$, which are non-overlapping, we obtain a tagged family in $T_\epsilon$, $\calP \defeq \bigcup_{j=1}^p \calP_j$, which is also a tagged family in $T$. Furthermore, $\calP$ is $\eta$-regular, $\delta$-fine and satisfies
\[
     G_i(T-[\calP]) \leq G_i(T-T_\epsilon) + \sum_{j=1}^q G_i(T_j-[\calP_j]) < \epsilon
\]
for $i=1,\dots,p$. Thus $T$ has the Cousin-Howard Property for the function $\eta_T$.
\end{proof}
In particular, we can apply this result to the case where $E_T$ is empty:
\begin{Corollary}\label{C1 has HC}
Currents associated to compact oriented $C^1$ submanifolds with boundary have the Cousin-Howard Property.
\end{Corollary}

\section{A generalized Stokes' Theorem and a counterexample}\label{sec:stokes}

We decompose Theorem \ref{thm:stokes} into two statements: the continuous case (Theorem \ref{thm:stokes2}) and the generalization to discontinuous forms (Proposition \ref{prop:stokesdiscont}). The last paragraph of this section contains the proof of Theorem \ref{thm:contrex}.
\subsection{Reduction to the continuous case}
\begin{Theorem}\label{thm:stokes2}
  Let $T\in \Imrn$ be weakly regular with singular set $E_T$, and $\omega: \spt T \to \Lambda^{m-1}(\Rn)$ be a differential form satisfying
  \begin{enumerate}[label=(\roman*)]
    \item $\omega$ is continuous on $\spt T$,\label{thm:stokes2:continuity}
    \item $\omega$ is pointwise Lipschitz continuous on $\spt T\backslash (E_\omega\cup E_T)$, where $E_\omega$ is $\scrH^{m-1}$ $\sigma$-finite,\label{thm:stokes2:lipschitz}
    \item $x\mapsto \langle \dd \omega(x), \vect T(x) \rangle$ is defined $\Vert T\Vert$ almost everywhere and Lebesgue integrable with respect to $\Vert T\Vert$.\label{thm:stokes2:integrability}
  \end{enumerate}
 Then there holds
\begin{equation}\label{eq:stokesstatement}
     \partial T (\omega)=\int \langle \dd \omega(x),\vect T (x) \rangle \dd \Vert T\Vert (x).
\end{equation}
\end{Theorem}
Assuming that the above holds, Theorem \ref{thm:stokes} follows from the following statement:
\begin{Proposition}\label{prop:stokesdiscont}
  Given the assumptions of Theorem \ref{thm:stokes2}, with \ref{thm:stokes2:continuity} replaced by
  \begin{enumerate}
    \item[(i')] $\omega$ is bounded on $\spt T$ and continuous on $\spt T\backslash E_0$, where $E_0$ is strongly disposable in $T$.
  \end{enumerate}
Then \eqref{eq:stokesstatement} holds.
\end{Proposition}
\begin{proof}
Let $M>0$ be an upper bound on $\vert \omega\vert$. Fix $\epsilon >0$, and use the strong disposability of $E_0$ to find $T_\epsilon \sqsubset T$ such that
\begin{equation*}
  \begin{cases}
    \spt T_\epsilon \cap E_0= \emptyset,\\
    \mass(\partial (T-T_\epsilon)) < \epsilon,\\
    \mass(T-T_\epsilon) < \epsilon.
  \end{cases}
\end{equation*}
There holds $\vert \partial T(\omega)-\partial T_\epsilon (\omega) \vert< M\epsilon$, so that when $\epsilon$ tends to zero, $\partial T_\epsilon(\omega)$ converges to the left hand side of \eqref{eq:stokesstatement}. Similarly, the integral of $\langle \dd \omega,\vect T\rangle$ with respect to $\Vert T_\epsilon\Vert$ is arbitrarily close to the right hand side of \eqref{eq:stokesstatement}. 

Clearly $\omega$ is continuous on $\spt T_\epsilon$, thus, supposing that $T_\epsilon$ is weakly regular, we can apply Theorem \ref{thm:stokes2} to $T_\epsilon$ and $\omega$. Letting $\epsilon$ go to zero yields \eqref{eq:stokesstatement}. The weak regularity of $T_\epsilon$ follows from Lemma \ref{lemma:weakreg-hered}, below.
\end{proof} 

It is not a priori clear that a subcurrent of a weakly regular current is weakly regular. Indeed, disposability is in general \textit{not} hereditary, in the sense that there exists a current $T$ and a set $E\subset \spt T$ which is disposable for $T$ but not for a certain subcurrent $S$ of $T$ (see Section 4.4.3 in \cite{ThesisAJ}). This example is linked to a current with small non disposable singular set similar to the one we define in the next section (Paragraph \ref{sec:example}). However, weak regularity is in fact hereditary:
\begin{Lemma}\label{lemma:weakreg-hered}
  If $T\in \Imrn$ is weakly regular with singular set $E_T$, then $E_T$ is disposable in every subcurrent $S$ of $T$.
\end{Lemma}
\begin{proof}
Let $T$ be a weakly regular current. Pick a countable locally finite family of charts $((\theta_j,A_j,\phi_j,U_j))_j$ of $T$ covering $\spt T\backslash E_T$ and let $(\psi_j)_j$ be a partition of unity subordinate to the $U_j$. Note that given a current $R$ supported in the union of the $U_j$, there holds $\partial R= \sum_j ((\partial R)\hel U_j)\hel \psi_j$.

Fix a subcurrent of $T$: $S\defeq  T\hel B$, and $\epsilon>0$. By the weak regularity of $T$, we can find a subcurrent $T_\epsilon$ of $T$ such that
\begin{equation*}
  \begin{cases}
    \spt T_\epsilon \cap E_T=\emptyset,\\
     \mass(T-T_\epsilon) <\epsilon,\\
     \mass(\partial T_\epsilon)<C,
  \end{cases}
\end{equation*}
where $C$ is independent of $\epsilon$. We can write $T_\epsilon = T\hel B_\epsilon$ for some set $B_\epsilon$. Consider the current $S_\epsilon \defeq S\hel B_\epsilon = T\hel (B\cap B_\epsilon)$, it satisfies $\mass (S-S\hel B_\epsilon) \leq \mass(T\hel B_\epsilon^c)<\epsilon$ and we claim that $\mass(\partial (S\hel B_\epsilon))\leq \mass (\partial T_\epsilon)+\mass(\partial S)$. As $E_T\cap B_\epsilon = \emptyset$, as $\epsilon$ is arbitrary, this would suffice to prove that $E_T$ is disposable in $S$.

There holds $\partial S_\epsilon= \sum_j (\partial S_\epsilon) \hel \psi_j$, so we only need to prove that for all $j$ 
\begin{equation}\label{eq:measurecomparison}
\Vert \partial S_\epsilon\Vert \hel U_j \leq \Vert \partial (T\hel B_\epsilon)\Vert\hel U_j + \Vert \partial S\Vert \hel U_j.
\end{equation}
Fix $j$ corresponding to a chart of $T$; there are charts of $T_\epsilon$ and $S$, such that we can write
\begin{eqnarray*}
  T_\epsilon\hel U_j &=& \theta_j {\phi_j} _\# \bE^m\hel  A^1_j,\\
 S\hel U_j &=& \theta_j {\phi_j}_\# \bE^m\hel A^2_j,
\end{eqnarray*}
where $A^1_j$ and $A_j^2$ are subsets of $A_j$ of finite perimeter (possibly of zero Lebesgue measure).
Noting that $A_j^1\cap A_j^2$ is a set of finite perimeter, we have $S_\epsilon \hel U_j = \theta_j {\phi_j}_\# \bE^m\hel( A^1_j\cap A_j^2)$, and we can also write
\begin{eqnarray*}
  (\partial T_\epsilon)\hel U_j &=& (\theta_j {\phi_j}_\# \partial (\bE^m\hel A_j^1))\hel U_j,\\
  (\partial S)\hel U_j &=& (\theta_j {\phi_j}_\# \partial (\bE^m\hel A_j^2))\hel U_j,\\
  (\partial S_\epsilon)\hel U_j &=&(\theta_j {\phi_j}_\# \partial(\bE^m \hel (A_j^1\cap A_j^2)))\hel U_j.
\end{eqnarray*}
By the standard theory of sets of finite perimeter (see for instance \cite[Chapter~2]{AmFuPa2000} or \cite[Chapter~5]{evans1991measure}), we have 
\[
     \Vert \partial (\bE^m\hel( A_j^1\cap A_j^2))\Vert \leq \Vert \partial (\bE^m \hel A_j^1)\Vert+ \Vert \partial (\bE^m\hel  A_j^2)\Vert,
\]
which directly implies inequality \eqref{eq:measurecomparison} and completes the proof.
\end{proof}
\subsection{Proof of the continuous case}
We now turn to the proof of Theorem \ref{thm:stokes2}. If $f$ is a function defined on $A\subset \spt T$ and $\calP$ is a tagged family in $T$ based in $A$, we will be interested in the \emph{Riemann sum of $f$ over $\calP$}:
\[
      \sigma(f,\calP) = \sum_{(x,S)\in\calP} f(x) \mass(S).
\]
Our proof of Stokes' Theorem is based on three building blocks. Lemma \ref{lemma:saks-henstock}, known as the Saks-Henstock Lemma allows us to approximate the Lebesgue integral of a function with respect to $\Vert T\Vert$, by appropriate Riemann sums associated to tagged families in $T$. The second (Lemma \ref{lemma:derivation}) is a derivation result, which specifies when the exterior differential of a form $\omega$ at a point $x$ approximates its circulation $\Theta_\omega$. The third result (Lemma \ref{lemma:negligibleset}) gives a way to control the circulation of $\omega$ near the points where it is non differentiable, but only pointwise Lipschitz.

\begin{Lemma}[Saks-Henstock Lemma]\label{lemma:saks-henstock}
  Suppose that $f:\spt T \to \R$ is Lebesgue-integrable with respect to $\Vert T\Vert$. Given $\epsilon>0$, there exists a positive gauge $\delta_1$ on $\spt T$ as well as a real number $\tau>0$ such that for any $\delta_1$-fine, $(\mass,\tau)$-full tagged family $\calP$ in $T$ there holds
  \begin{equation}
    \label{eq:saks-henstock}
    \left \vert \int f \dd \Vert T\Vert - \sigma(f,\calP) \right \vert < \epsilon.
  \end{equation}
\end{Lemma}

\begin{proof}
 We first fix a representative of $f$ on $\spt T$. As the measure $\Vert T\Vert$ is finite and Borel regular we can apply the Vitali-Caratheodory Theorem to $f$ (see \cite[2.24]{RudinReal}).  Choose $\epsilon >0$. There exist extended-real valued functions $g$ and $h$ defined on $\spt T$, which are respectively upper and lower semi-continuous, such that $g\leq f\leq h$  and satisfy 
\begin{align*}
\int (h-f)\dd \Vert T\Vert <\epsilon/2\quad \text{ and } \quad\int (f-g)\dd \Vert T\Vert<\epsilon/2.
\end{align*}
 Now, using the upper semi-continuity of $g$ and the lower semi-continuity of $h$, for each $x \in \spt T$, we fix a positive $\delta_1(x)$ such that for all $y\in \cBall(x,\delta(x))$
\begin{equation}\label{eq:SHcontrol1}
    g(y) -\dfrac{\epsilon}{4\mass(T)} \leq g(x) \leq f(x) \leq h(x)\leq h(y) +\dfrac{\epsilon}{4\mass(T)}.
\end{equation}
  Let $S$ be a subcurrent of $T$ and suppose that for some  $x\in \spt S$, there holds $\diam (\spt S) <\delta_1(x)$, then there holds
\[
    \int g \dd \Vert S\Vert -\dfrac{\epsilon\mass(S)}{4 \mass(T)} \leq  f(x)\mass(S) \leq \int h \dd \Vert S \Vert + \dfrac{\epsilon \mass(S)}{4 \mass(T)}.
\]
Thus if $\calP$ is a $\delta_1$-fine tagged family in $T$, we have
\begin{equation}\label{eq:SHcontrol2}
\int g \dd \left  \Vert [\calP]\right \Vert-\dfrac{\epsilon}{4}\leq \sigma(f,\calP) \leq \int h \dd \left  \Vert [\calP]\right \Vert + \dfrac{\epsilon}{4}.
\end{equation}
By the Lebesgue integrability of $g$ and $h$, we can find $\tau>0$ such that if $R$ is a subcurrent of $T$ with $\mass(T-R) <\tau$, then
\begin{equation}\label{eq:SHcontrol3}
     \int \vert g\vert \dd \Vert T-R\Vert<\dfrac{\epsilon}{4},
\end{equation}
and the same estimate holds when integrating $\vert h\vert$. Thus if $\calP$ is a $\delta_1$-fine tagged family in $T$ with $\mass(T-[\calP])<\tau$, estimates \eqref{eq:SHcontrol2} and \eqref{eq:SHcontrol3} yield
\begin{equation*}
\int g \dd  \Vert T \Vert-\dfrac{\epsilon}{2}\leq \sigma(f,\calP) \leq \int h \dd  \Vert T \Vert + \dfrac{\epsilon}{2},
\end{equation*}
 which combined with  \eqref{eq:SHcontrol1} implies \eqref{eq:saks-henstock} and concludes the proof of the Saks-Henstock Lemma.
\end{proof}
\begin{Lemma}[Derivation Lemma]\label{lemma:derivation}
Fix $T\in \Imrn$, suppose that $x\in \set_m\Vert T\Vert$  is a regular point of $T$ and that $\omega:\spt T\to \Lambda^{m-1}(\Rn)$ is differentiable at $x$. Then for all $\eta>0$ and $\epsilon>0$, there exists a real number $\delta_2(x)>0$ such that whenever $S$ is an $\eta$-regular subcurrent of $T$ with $x \in \spt S$ and $\diam \spt S <\delta_2(x)$, there holds
  \begin{equation}
    \label{eq:derivation}
    \left \vert \langle \dd \omega(x), \vect T(x)\rangle \mass(S) - \int \langle \omega,\vect{\partial S} \rangle \dd \Vert \partial S \Vert \right \vert <\epsilon \mass (S).
  \end{equation}
\end{Lemma}
\begin{Remark}\label{rem:appcont}
    I think that the continuity assumption on $\vect T$ at $x$ is sharp if $1<m<n$. Indeed, the example for \cite[Conjecture 3.2.17]{ThesisAJ} tends to show that the approximate continuity of $\vect T$ is not sufficient. However, if in a neighbourhood of $x$, $T$ is the bi-Lipschitz push-forward of some set of finite perimeter, approximate continuity of $\vect T$ might be sufficient for \eqref{eq:derivation} to hold, but I was not able to prove it.
  \end{Remark}
  \begin{proof}[Proof of Lemma \ref{lemma:derivation}]
Fix $\eta>0$ and $\epsilon>0$. As $\omega$ is differentiable at $x$, there exists $\delta>0$ such that for $y\in \spt T\cap \Ball(x,\delta)$, there holds
\begin{equation}\label{eq:diffomega}
   \vert \omega(y)-\omega(x) - (y-x)\leh \dd \omega (x) \vert <\dfrac{\epsilon}{2\eta} \vert y-x\vert,
 \end{equation}
where $\leh$ represents the inner product in the notation of \cite{FedererGMT}. Furthermore since $x$ is a regular point of $T$, $\vect T$ has a continuous representative at $x$ on $\spt T$, thus we can suppose that for $y\in \spt T\cap \Ball(x,\delta)$
\begin{equation}\label{eq:deriv3}
    \vert \vect{T}(y)-\vect{T}(x)\vert <\dfrac{\epsilon}{2 \max\{\vert \dd \omega(x)\vert, 1\}}.
  \end{equation}
  Let $S$ be an $\eta$ regular subcurrent of $T$, with $x\in \spt S$ and $\diam(\spt S)<\delta$. A consequence of \eqref{eq:deriv3} is that $S$ is not a cycle. To see this, test $\partial S$ against a restriction of the $m-1$ form $\zeta: y\mapsto (y-x)\leh (\vect{T}(x))^*$ (where if $v$ is a $k$-vector, $v^*$ stands for its dual $k$-covector). There holds $\dd \zeta (y)= (\vect{T}(x))^*$ for all $y\in \Rn$ and it is easy to show that $\partial S(\zeta) = S(\dd \zeta) \neq 0$.

  By \eqref{eq:diffomega}, we can approximate the integral of $\omega$ on $\partial S$ as follows:
\begin{align}\label{eq:deriv1}
   \left \vert \int \left ( \langle \omega(y),\vect{\partial S}(y) \rangle - \langle \omega(x)  (y-x)\leh \dd \omega (x),\vect{\partial S}(y)\rangle \right )\dd \Vert \partial S\Vert (y)\right \vert \nonumber\\ < \dfrac{\epsilon \diam (\spt S)\mass(\partial S)}{2\eta} < \dfrac{\epsilon\mass( S)}{2}.
\end{align}
As the $(m-1)$-form $y\mapsto \omega(x) + (y-x) \leh \dd \omega(x)$ is smooth and has constant differential equal to $\dd \omega(x)$, the Stokes' theorem for smooth forms and integral currents implies
\begin{equation}\label{eq:deriv2}
  \int \langle \omega(x) + (y-x)\leh \dd \omega (x),\vect{\partial S}(y)\rangle \dd \Vert \partial S\Vert = \int \langle \dd \omega(x), \vect T(y) \rangle\Vert S\Vert (y).
\end{equation}
Using estimate \eqref{eq:deriv3} in \eqref{eq:deriv2}, we get
\begin{equation}
  \label{eq:deriv4}
  \left \vert \int \langle \omega(x) + (y-x)\leh \dd \omega (x),\vect{\partial S}(y)\rangle \dd \Vert \partial S\Vert- \langle \dd \omega(x), \vect T(x) \rangle \mass(S)\right \vert < \dfrac{\epsilon}{2}\mass(S).
\end{equation}
Finally estimate \eqref{eq:derivation} follows for $\delta_2(x)\defeq \delta$ by combining \eqref{eq:deriv1} and \eqref{eq:deriv4}.
\end{proof}
We now treat the points where $\omega$ is not differentiable.
\begin{Lemma}\label{lemma:negligibleset}
  Let $\omega$ be a differential form of degree $m-1$, which is continuous on $\spt T$ and pointwise Lipschitz on $\spt T \backslash (E_T\cup E_\omega)$, and let $\eta>0$ be a positive function on $\spt T\backslash (E_T\cup E_\omega)$; there holds
  \begin{enumerate}
    \item $\omega$ is differentiable on $\spt T\backslash (N\cup E_\sigma \cup E_T) $, where  $N\cap ( E_\sigma\cup E_\omega) = \emptyset$ and $\Vert T\Vert(N)= 0$.\label{item:aediff}
    \item\label{item:negligible} Given $\epsilon>0$, there exists a positive function $\delta_3$ on $N$ such that if $\calP$ is a $\delta_3$-fine and $\eta$-regular tagged family, then 
\[   
\vert (\partial [\,\calP\,])(\omega) \vert < \epsilon \quad \text{ and } \quad\mass([\,\calP\,]) < \epsilon.
\]
 \end{enumerate}
\end{Lemma}
\begin{proof}
  To prove \eqref{item:aediff}, it suffices to work in a $C^1$ chart, as $\spt T\backslash E_T$ can be covered by countably many charts, and a countable union of $\Vert T\Vert$-null sets is also $\Vert T\Vert$-null. Let $U$ be an open set corresponding to a chart of $T$, such that $\spt T\cap U$ is a $C^1$ submanifold of dimension $m$ in $\Rn$. In $(\spt T\cap U)\backslash E_\omega$, $\omega$ is pointwise Lipschitz and is thus $\omega$ is differentiable $\scrH^m\hel(\spt T\cap U \backslash E_\omega)$ almost everywhere by the Rademacher-Stepanov Theorem (see \cite[3.1.6 and 3.1.9]{FedererGMT}). As $\Vert T\Vert\hel U$ is absolutely continuous with respect to $\scrH^m\hel (\spt T\cap U)$, the claim is proved. Let $N$ be the set of non-differentiability points of $\omega$ in $\spt T\backslash (E_\omega \cup E_T)$.

Let us now prove claim \eqref{item:negligible}. Fix $\epsilon>0$ and notice that given $x\in \spt T\backslash ( E_T\cup E_\omega)$, there exists $\delta(x)>0$ such that for $S\sqsubset T$ with $\diam \spt S < \delta(x)$, $x\in \spt S$ and $\mass(S) >\eta(x) \mass(\partial S) \diam (\spt S)$, there holds
\begin{equation}\label{eq:controlLip}
       \vert(\partial S)(\omega ) \vert \leq \Lip_x \omega\dfrac{\mass(S)}{\eta(x)}.
\end{equation}
     For $k=1,2,\dots$, let $N_k\defeq \{x\in N, k-1\leq  \Lip_x \omega / \eta(x) <k\}$. As $\Vert T\Vert$ is a Radon measure, we can find an open set $U_k$ containing $N_k$ and with $\Vert T\Vert ( U_k ) < k^{-1} 2^{-k}\epsilon$. Taking $\delta(x)$ smaller if necessary, for each $x\in U_k$ we can suppose that if $S\sqsubset T$ with $x\in \spt S$ and $\diam \spt S <\delta(x)$, then $\spt S \subset U_k$. Fix $\delta_3= \delta$, it is a positive function on $N$. There remains to be shown that the conclusions of \eqref{item:negligible} hold. Let $\calP$ be a $\delta_3$-fine, $\eta$-regular family in $T$, we have
\[
      \mass([\,\calP\,]) \leq \sum_{k= 1}^{\infty} \Vert T\Vert (U_k) < \epsilon.
\]
and
\[
    \vert (\partial[\, \calP\,]) (\omega)\vert \leq \sum_{(x,S)\in \calP} \vert(\partial S)(\omega )\vert \leq \sum_{k=1}^\infty\sum_{(x,S)\in \calP,\, x\in N_k}  k \mass(S ) \leq \sum_{k=1}^{\infty} k\Vert T\Vert(U_k) < \epsilon,
\]
 by estimate \eqref{eq:controlLip} and the definition of $N_k$.
\end{proof}

\begin{proof}[Proof of Theorem \ref{thm:stokes2}]
  Let us first recall, that as seen in Example \ref{example:additive}, the function $\Theta_\omega: S\mapsto \partial S(\omega)$ for $S\in \subcspace (T)$ is continuous and additive. Additionally, by Theorem \ref{thm:weaklyreg-hc},  $T$ has the Cousin-Howard Property for some positive function $\eta$. Furthermore, by Lemma~\ref{lemma:negligibleset}~\eqref{item:aediff}, $\omega$ is differentiable along $\spt T$, $\Vert T\Vert$ almost everywhere, we let $N$ be the set of points of $\spt T\backslash (E_\sigma\cup E_T)$ where $\omega$ is not differentiable.

  Fix $\epsilon>0$. Consider the positive function $\delta_1$ and the positive number $\tau$ associated to  the integrable function $\langle \dd \omega,\vect T\rangle$ and this choice of $\epsilon$ by Lemma \ref{lemma:saks-henstock}. Consider also the function $\delta_2$ defined at a regular point $x$ of $T$ at which $\omega$ is differentiable, by using Lemma \ref{lemma:derivation} with the parameters $\eta(x)$ and $\epsilon$. Lastly consider the function $\delta_3$ on $N$, corresponding to $\min\{\epsilon,\tau/2\}$ and $\eta$ in Lemma~\ref{lemma:negligibleset}~\eqref{item:negligible}. We define a gauge $\delta$ on $\spt T$ as follows:
  \begin{equation*}
    \delta(x) = \begin{cases}
      \min\{\delta_1(x),\delta_2(x)\} &\text{ if } x \in \spt T \backslash (E_\sigma \cup E_T\cup N),\\
      \delta_3(x) &\text{ if } x\in N,\\
      0 &\text{ if } x\in E_\sigma \cup E_T.
      \end{cases}
  \end{equation*}
Note that $\delta$ is positive on $\set_m T$ except on a $\scrH^{m-1}$ $\sigma$-finite set. We can thus apply the Cousin-Howard Property to $T$ for the continuous additive functions $\Theta_\omega$, $-\Theta_\omega$ and $\mass\vert_{\subcspace(T)}$, with the parameter $\min\{\epsilon, \tau/2\}$. This yields a $\delta$-fine, $\eta$ regular, tagged family $\calP$ in $T$ such that
\begin{equation}\label{eq:circapprox}
      \begin{cases}
        \mass(T-[\calP]) &<\tau/2,\\
        \vert \Theta_\omega(T-[\calP])\vert &<\epsilon.
      \end{cases}
    \end{equation}
Now, $\calP$ can be partitionned into $\calP_N\cup \calP_d$, where $\calP_N$ consists of all pairs $(x,S)\in \calP$ such that $x\in N$. By the choice of $\delta_3$ and Lemma \ref{lemma:negligibleset}, there holds 
\begin{equation}\label{eq:neglig}
  \begin{cases}
     \mass([\,\calP_N\,])<\tau/2\\
\vert \Theta_\omega([\,\calP_N\,])\vert <\epsilon.
  \end{cases}
\end{equation}
This implies in particular that $\mass(T-[\,\calP_d\,])<\tau,$ so applying Lemma \ref{lemma:saks-henstock} to $\calP_d$ and the function $f= \langle \dd \omega,\vect T\rangle$ yields 
\begin{equation}\label{eq:riemannsum1}
    \left \vert \int \langle \dd \omega,\vect T\rangle \dd \Vert T\Vert - \sigma(f,\calP_d) \right \vert < \epsilon.
\end{equation}
Since each $(x,S)\in \calP_d$ satisfies the conditions of Lemma \ref{lemma:derivation}, estimate \eqref{eq:derivation} holds, and summing over the elements of $\calP_d$ yields
\begin{equation}
    \label{eq:riemannsum2}
    \left \vert \Theta_\omega([\calP_d]) - \sigma(f,\calP_d)\right \vert <\epsilon \mass ([\calP_d]).
  \end{equation}
Combining estimates \eqref{eq:circapprox}, \eqref{eq:neglig}, \eqref{eq:riemannsum1} and  \eqref{eq:riemannsum2} yields
\begin{equation*}
    \left \vert \int \langle \dd \omega,\vect T\rangle \dd \Vert T\Vert -  \Theta_\omega(T) \right \vert <(2 + \mass(T)) \epsilon.
\end{equation*}
As $\epsilon$ is arbitrary, this completes the proof of Theorem \ref{thm:stokes2}.
\end{proof}

\subsection{A current which is not weakly regular}\label{sec:example}
In this section, we provide examples of currents whose singular sets have small Hausdorff measure and dimension, yet possibly infinite Minkowski content. This alone does not mean that their singular sets are disposable. However we also give an example of current of dimension $2$ whose singular set is reduced to a point, but not disposable. This will prove the following result (already stated in the introduction):

\contrex*

For clarity, we outline the construction first, and postpone the technical details until the end of the proof.
\begin{proof}[Proof of Theorem \ref{thm:contrex}]
  We construct the current $T$ as the graph in $\R^3$ over the disk $D\defeq \Ball_{\R^2}(0,3/2)$ of a function $f$ defined on concentric annuli as follows:

For $k=0, 1,2, \dots$, let $r_k \defeq \sum_{j= k}^\infty 3^{-j}$ (in particular $r_0= 3/2$). The crown situated between the radii $r_{k+1}$ and $r_k$ has thickness $3^{-k}$. In polar coordinates, we define $f_0(\theta)= 0$ and for $k\geq 1$, $f_k(\theta) = 3^{-k} \sin (4^{k}\theta)$. Let $\varphi: [\,0,1\,]\to [\,0,1\,]$ be a smooth function, which is equal to $1$ on $[\,0,1/4\,]$ and to $0$ on $[\,3/4,1\,]$ and satisfies $\vert \varphi'\vert\leq 4$. For $r\in [\,r_{k+1},r_k\,]$ and $\theta \in [\,0,2\pi)$, we let 
\[
     f(r,\theta) \defeq \varphi (3^{k}(r_k-r) ) f_k(\theta) + \varphi \left (3^k(r-r_{k+1})\right ) f_{k+1}(\theta).
\]
We do this for all $k$, in order to define $f$ on $D\backslash \{0\}$. Letting $f(0,0)=0$ we extend $f$ by continuity to the whole of $D$. Let $F:D\to \R^3$ be the graph map given in cylindrical coordinates by $F(r,\theta) = (r,\theta,f(r,\theta))$. Let $M$ be the image of $F$.

The function $f$ is smooth outside of the origin and $M\backslash\{0\}$ is an oriented submanifold of dimension $2$ in $\R^3$. We can prove that $M$ has finite area. Thus $M$ can be represented by a rectifiable current $T$ of dimension $2$. We can also show that $T$ is integral.
\begin{Claim}\label{claim:defineT}
  The {\it push-forward} by $(r,\theta)\mapsto (r,\theta,f(r,\theta))$ of $\bE_2 \hel \cBall_{\R^2}(0,3/2)$, with a slight abuse of notation as the map is not Lipschitz continuous at the origin, defines an integral current $T$ whose boundary represents a circle. The only singular point of $T$ is the origin.
\end{Claim}
The following important feature of $M$ is the key to our example:
\begin{Claim}\label{claim:length}
  The length $L(R)$ of the curve $M\cap \{(r,\theta,z), r=R\}$ goes to infinity as $R$ goes to zero.
\end{Claim}
In particular $\mass (\partial (T\hel \{r<R\})) \to \infty$ as $R$ goes to zero. However, this does not imply that the singular set of $T$ is not disposable. 

We prove that $E_T = \{0\}$ is not disposable by contradiction with Theorem \ref{thm:stokes}. More precisely, we define a continuous form $\omega: \spt T\to \Lambda^1(\R^3)$, which is differentiable at all points of $\spt T\backslash \{0\}$ and whose differential along $T$ is zero. In particular, this differential is integrable with respect to $\Vert T\Vert$. However we construct $\omega$ in such a way that $T(\omega)\neq 0$ contradicting \eqref{eq:stokes}.

This form is constructed as the pullback by the orthogonal projection $\pi:M\to D$ of the differential of the smooth map $u: D\backslash\{0\}\to \R/\Z$ defined by
\[
     u(r,\theta) = \left (\dfrac{L(r,\theta)}{L(r)}\right )/\{0\sim 1\},
\]
where $L(r,\theta)$ is the length of the portion of curve $M\cap \{(r,\theta',z), \theta' \in (0,\theta)\}$.
\begin{Claim}\label{claim:extendomega}
  The form $\omega$ is well defined on $M\backslash \{0\}$, it can be extended by continuity to $0$ at $0$, and satisfies $\langle \dd \omega(x),\vect T(x) \rangle=0$ for $x\in M\backslash\{0\}$.
\end{Claim}
To conclude, notice that $\partial T(\omega) = L(3/2)/L(3/2)-0 = 1.$ This implies that \eqref{eq:stokes} does not hold, and therefore $\{0\}$ is not disposable in $T$. 
\end{proof}
We now outline the proofs of the three claims above
\begin{proof}[Proof of Claim \ref{claim:defineT}]
  The only difficult point is to prove that $T$ has finite mass, or equivalently, that $\scrH^2(M)<+\infty$. Let $M_k$ be the image by $F$ of the crown situated between the radii $r_{k+1}$ and $r_k$, there holds:
  \begin{align*}
    \scrH^2(M_k)&= \int_0^{2\pi} \int_{r_{k+1}}^{r_k} \sqrt{1+ ( \partial_r f(r,\theta))^2+ \left(\dfrac{\partial_\theta f(r,\theta)}{r}\right)^2 }r\dd r \dd \theta\\
    &\leq 2\pi (r_k-r_{k+1}) r_k \sqrt{1+ (3^{-k} 3^k \max (\vert \varphi'\vert))^2+( r_k^{-1}(3^{-k}4^k+ 3^{-k-1}4^{k+1}))^2}\\
    &\leq C 3^{-2k}\sqrt{1+ 4^{2k}}.
  \end{align*}
In particular $\sum_k \scrH^{2}(M_k) <+\infty$.

We infer that $T$ is a rectifiable current of dimension $2$ (with finite mass by definition), such that $\partial T$ is supported on $F(\partial D)\cup \{0\}$. As a flat chain of dimension $1$ cannot be supported on a set of dimension $0$, $\partial T$ is supported on $F(\partial D)$, which is a connected closed curve, and $T$ is thus integral.
\end{proof}
\begin{proof}[Proof of Claim \ref{claim:length}]
  The length of the curve $M\cap \{(r,\theta,z), r=R\}$ is given by
    \begin{align*}
      L(R) = \int_{0}^{2\pi} \sqrt{ 1+ \left(\dfrac{\partial_\theta f(R,\theta)}{R}\right)^2}R\dd \theta
      \geq 2 \pi \int_{0}^{2\pi} \vert\partial_\theta f(R,\theta)\vert\dd \theta.
    \end{align*}
If $R = r_{k+1} + 3^{-k}s$ with $s\in [\,0,1\,]$ ($s= 1$ yields $R= r_k$), we can write
\[
    \partial_\theta f(R,\theta) = \varphi(1-s)3^{-k} 4^k\cos(4^{k}\theta) + \varphi(s) 3^{-k-1}4^{k+1}\cos(4^{k+1}\theta).
\]
Integrating the absolute value over $\theta$, yields
\[
     L(R) \geq 4^{k}3^{-k}g(s),
\]
where for some positive continuous function $g$ defined on $[\,0,1\,]$ and depending only on $\varphi$. As $k\to \infty$ when $R\to 0$, we have $L(R)\to +\infty$.
\end{proof}
\begin{proof}[Proof of Claim \ref{claim:extendomega}]
It is clear that $\omega$ is locally closed in $M\backslash \{0\}$, as the pullback of an exact form. We need only to prove that $\omega$ can be extended by continuity at the origin.
Notice that $\dd u$ cannot be bounded, for otherwise it would violate Theorem \ref{thm:stokes} on $\bE_2\hel D\in \Int_2(\R^2)$, which is weakly regular. So we have to estimate $\langle \omega, \bv\rangle $ for $\bv$ tangent to $M$.

 For such a vector, we have $\langle \omega, \bv \rangle = \langle \dd u, \pi (\bv)\rangle$.
First, let us express $\dd u$. For $(r,\theta) \in D\backslash \{0\}$, supposing for simplicity that $\theta \neq 0 \mod 2\pi$, we have
  \begin{align*}
      \dd u(r,\theta) &= \partial_r u (r,\theta) \be^*_r + r^{-1} \partial_\theta u(r,\theta) \be^*_\theta \\
&=  \dfrac{L(r) \partial_r L(r,\theta) - L(r,\theta) \partial_r L(r)}{L(r)^2} \be^*_r +  \dfrac{2\pi\sqrt{1+r^{-2}(\partial_\theta f)^2(r,\theta)}}{L(r)} \be_\theta^*. 
  \end{align*}
    Let us choose a tangent orthonormal basis to $M$ at the point $(r,\theta,f(r,\theta))$, we proceed by orthonormalization of the basis of $T M$: $(\dD F\,\be_\theta, \dD F\, \be_r)$. In the basis $(\be_r,\be_\theta,\be_z)$, this yields
\begin{align*}
   \tau_1(r,\theta)&= \dfrac{1}{\sqrt{1+r^{-2}(\partial_\theta f)^2}}\begin{pmatrix}
    0\\1\\ r^{-1}\partial_\theta f \end{pmatrix},\\
  \tau_2 (r,\theta) &= \dfrac{1}{\sqrt{1+r^{-2}(\partial_\theta f)^2}\sqrt{1+r^{-2}(\partial_\theta f)^2+(\partial_r f)^2}}\begin{pmatrix}
     1+r^{-2}(\partial_\theta f)^2\\ -r^{-1} \partial_\theta f \partial_r f\\ \partial_r f
    \end{pmatrix}.
\end{align*}
There holds
\begin{align*}
  \langle \dd u, \pi(\tau_1) \rangle = \dfrac{2\pi}{L(r)},
\end{align*}
thus $\vert \langle \omega(r,\theta,f(r,\theta) ,\tau_1(r,\theta) \rangle \vert \leq 2\pi L(r)^{-1}\to 0$, as $r$ goes to $0$. 

The term in $\tau_2$ is slightly harder to control, we have
\begin{align*}
  \langle \dd u, \pi(\tau_2) \rangle &= \dfrac{-2\pi (r^{-1} \partial_\theta f \partial_r f) }{L(r)\sqrt{1+r^{-2}(\partial_\theta f)^2+(\partial_r f)^2}}\\
  &\hspace{2cm} +\dfrac{L(r) \partial_r L(r,\theta) - L(r,\theta) \partial_r L(r)}{L(r)^2} \dfrac{\sqrt{1+r^{-2}(\partial_\theta f)^2}}{\sqrt{1+r^{-2}(\partial_\theta f)^2+(\partial_r f)^2}}.
\end{align*}
As $\vert \partial_r f\vert \leq 2 \max\vert \varphi '\vert$, the first member of the sum is clearly controlled by $L(r)^{-1}$. To control the second term, it suffices to control its first factor which is equal to $\langle \dd u,\be_r^*\rangle$. We do so using the periodicity of $f$ in the variable $\theta$.

Suppose that $r\in [\,r_{k+1},r_k\,]$ and that $\theta \in [\, j(2\pi 4^{-k}), (j+1)(2\pi 4^{-k})\,]$ for some integer $j<4^k$. Then we can write
\[
    L(r,\theta) = 4^{-k}j L(r) + L(r, \theta-2\pi 4^{-k} j),
\]
thus 
$$
\partial_r L(r,\theta) = 4^{-k} j \partial_r L(r) + \partial_r L(r,\theta-2\pi 4^{-k} j)
$$
 and we have
\begin{equation}\label{eq:contrexdiff}
  L(r) \partial_r L(r,\theta) - L(r,\theta) \partial_r L(r)  =  \partial_r L(r,\theta-2\pi 4^{-k} j) L(r) - \partial_r L(r) L(r,\theta-2\pi 4^{-k} j)
 \end{equation}
For $\alpha\in (0,2\pi\,]$, there holds
\begin{align*}
     \vert \partial_r L(r,\alpha)\vert  \leq \int_0^\alpha \dfrac{ r^{-2} \vert \partial_\theta f\vert \vert \partial^2_{r,\theta} f\vert }{\sqrt{1+r^{-2}(\partial_\theta f)^2}} r\dd \theta \leq \alpha  \max\{\vert \partial^2_{r,\theta} f(r,\cdot) \vert\} \leq C \alpha 4^k.
\end{align*}
Plugging the last estimate into \eqref{eq:contrexdiff} with $\alpha = 2\pi$ and $\alpha = \theta-2\pi 4^{-k}j < 2\pi 4^{-k}$, using the easy estimate $L(r,\theta-2\pi 4^{-k} j) \leq 4^{-k} L(r)$ and dividing by $L(r)^2$ yields 
\begin{align*}
  \vert \langle \dd u,\be_r^*\rangle\vert \leq \dfrac{C}{L(r)} \to 0 \quad \text{as} \quad r\to 0.
\end{align*}
This completes the proof that $\langle \omega,\tau_2\rangle$ goes to zero as $r$ goes to $0$. We infer that $\omega$ can be extended by continuity at the origin.
\end{proof}

\section{Minkowski content and weakly regular currents}\label{sec:Minkowski}

In this section we prove Proposition \ref{prop:minkowskidispo} and parts \ref{goodmink} and \ref{goodmini} of Theorem \ref{thm:goodcurrents}.

\subsection{Minkowski content and disposability}
\begin{Definition}\label{def:minkowski}
Given a current $T\in \Imrn$ and a set $E\subset \spt T$, the \emph{lower $m-1$ dimensional $\Vert T\Vert$-Minkowski content of $E$} is the extended real number
\begin{equation*}
    \mathscr{M}^{m-1}_{\Vert T\Vert,*} (E) \defeq \liminf_{r\to 0} \dfrac{\Vert T\Vert(\Ball(E,r))}{r}.
\end{equation*}
\end{Definition}

We first compare the lower $m-1$ dimensional $\Vert T\Vert$-Minkowski content with the Hausdorff measure of dimension $m-1$ on the density set of $\Vert T\Vert$:
\begin{Proposition}\label{prop:minkowski-comparison}
  There exists a positive constant $c$ depending only on $m$ and $n$ such that for $E \subset \spt T$, there holds
\[
   \mathscr{M}^{m-1}_{\Vert T\Vert,*}(E) \geq c \scrH^{m-1}(E\cap \set_m\Vert T\Vert).
\]
\end{Proposition}

\begin{proof}
  Writing $T= \theta \scrH^m\hel M \wedge \vect T$, we can suppose that $M= \set_m\Vert T\Vert$ and there holds for $r>0$:
\begin{equation}\label{eq:measureestim}
     \Vert T\Vert (\Ball(E,r)) \geq \scrH^m (\Ball(E,r)\cap M).
\end{equation}
   As in the usual comparison of Minkowski content and Hausdorff measure (see P. Mattila's book \cite[Chapter 5.3-5.5]{Mattila1995book}), we introduce the covering and packing numbers of a set $A\subset \Rn$: respectively $\calN(A,r)$ and $\calP(A,r)$ and recall that $\calP(A,r) \geq \calN(A,2r)$. Estimate \eqref{eq:measureestim} implies 
\begin{equation*}
    \Vert T\Vert (\Ball(E,r)) \geq \calP(E\cap M,r) \alpha_m r^m \geq \calN(E\cap M,2r) \alpha_m r^m.
\end{equation*}
  Given $s>0$, if $\scrH^{m-1}_s$ is the size $s$ approximating Hausdorff measure of dimension $m-1$, we have for some constant $c>0$ depending only on $m$ and $n$:
\[
      \Vert T\Vert (\Ball(E,r)) \geq\calN(E\cap M,2r)\alpha_m r^{m} \geq c r \scrH^{m-1}_{2r}(E\cap M).
\]
Therefore, dividing by $r$, letting $r$ go to zero and taking the lower limit of the left hand side, we get
\[
    \mathscr{M}^{m-1}_{\Vert T\Vert,*} (E) \geq c \lim_{r\to 0} \scrH^{m-1}_{2r} (E\cap M) = c\scrH^{m-1}(E\cap M).
\]
\end{proof}
A reverse inequality cannot hold, indeed, here are two examples showing that the relative Minkowski content can be much larger than the Hausdorff measure:
\begin{Example}\label{example:infmink}  
The current $T_1$ in $\R^2$ representing the graph the function $x\mapsto x^{3/2}\cos(x^{-1})$ for $x \in [\,-\pi^{-1},\pi^{-1}\,]$ is in $\Int_1(\R^2)$. Its only singular point is the origin. The singular set $\{(0,0)\}$ has infinite $0$-dimensional lower $\Vert T_1\Vert$ Minkowski content. However, the origin is clearly disposable in $T_1$.
\end{Example}
\begin{Example}
   Let $T_2$ be the current in $\Int_1 (\R^2)$ representing the union of the oriented circles with center $(2^{-k},0)$ and radius $2^{-k}$ for $k\geq 1$. Then $T_2$ is in $\Int_1(\R^2)$ and its only singular point is $(0,0)$. Furthermore, $\mathscr{M}^{0}_{\Vert T_2\Vert,*}(\{(0,0)\})= +\infty$. Clearly the origin is not disposable in $T_2$, however a generalized Stokes' Theorem holds on $T_2$, as is the case on all integral currents of dimension $1$ in Euclidean spaces. This is proved in \cite{julia1daccepted}, using the decomposition of integral currents of dimension $1$ into countable sums of simple Lipschitz curves.
\end{Example}
We can now prove the result which motivates the introduction of the intrinsic Minkowski content:
\minkowskidispo*
Note in particular that statement \eqref{minkdispo} implies directly part \ref{goodmink} of Theorem \ref{thm:goodcurrents}.
\begin{proof}[Proof of Proposition \ref{prop:minkowskidispo}]
We start with \eqref{minkdispo}. Fix $\epsilon>0$ and let $C\defeq\mathscr{M}^{m-1}_{\Vert T\Vert,*} (E)$. There exists $r_0>0$ such that for all $r\in (0,r_0\,]$,
\[
    \Vert T\Vert (\Ball(E,r))<  2C r_0 <\epsilon.
\]
Since $T$ is an integral current, for $r$ in a set of positive measure in $(r_0/2,r_0)$, the slice $\langle T,\dist(E,\cdot),r\rangle$ is an integral current and
\begin{align*}
    \mass ( \langle T,\dist(E,\cdot), r\rangle) \leq \dfrac{2}{r_0} \Vert T\Vert (\Ball(E,r_0))\leq 4C.
\end{align*}
From this we infer
\begin{align*}
    \mass(\partial(T\hel (\Ball(E, r)))) &\leq \Vert \partial T \Vert(\Ball(E,r))+  \mass(<T,\dist(E,\cdot),r>)\\
&\leq \mass(\partial T) + 4C. 
\end{align*}
Pick such and $r \eqdef r_1$ and let $T_\epsilon \defeq T\hel (\Ball(E,r_1)^c)$. There holds $\spt T_\epsilon\cap E = \emptyset$ and
\begin{align*}
  \mass(T-T_\epsilon)  &< \epsilon,\\
  \mass(\partial T_\epsilon) &\leq 4C+\mass(\partial T).
\end{align*}
Such a $T_\epsilon$ can be chosen for any $\epsilon>0$,  thus $E$ is disposable in $T$.

The proof of \eqref{minkstrongdispo} follows from a similar argument, using the additional fact that since $\mass(\partial T) <+\infty$ and $\Vert \partial T\Vert(E)=0$, there holds 
\[
    \Vert \partial T\Vert (\Ball(E,r))\to 0 \quad \text{ as } r\to 0.
\]
\end{proof}

\subsection{Mass minimizing integral currents}
In this paragraph, we prove statement \ref{goodmini} of Theorem \ref{thm:goodcurrents}, which in light of the previous paragraph, can be rephrased as:
\begin{Theorem}\label{thm:minimizing-weaklyreg}
  Suppose that $T$ is a mass minimizing current of dimension $n-1$ in $\Rn$ and that $\partial T$ represents a closed oriented $C^{1,\alpha}$ submanifold of $\Rn$ for some $\alpha>0$ with multiplicity one. Then $T$ is weakly regular.
\end{Theorem}

In light of part \ref{minkdispo} of Proposition \ref{prop:minkowskidispo}, we only need to prove that the singular set $E_T$ of $T$ has finite lower Minkowski content of dimension $n-2$ with respect to $\Vert T\Vert$. To do this we consider separately the interior and boundary parts of $E_T$. More precisely, let $E_T^\partial \defeq E_T\cap \spt \partial T$ and $E_T^i \defeq E_T\backslash E_T^\partial$. We start by obtaining a local estimate on $E_T^i$; recall that the classical regularity theory (\cite{BomGioGiu1969}) implies that $E_T^i$ has Hausdorff dimension at most $n-8$, however, this is not sufficient to infer a Minkowski content estimate.

A recent improvement to the interior regularity theory was obtained by J. Cheeger and A. Naber:
\begin{Theorem}[{\cite[Theorem 5.8]{CheegerNaber2013}}]\label{thm:CheegerNaber}
Suppose $T\in \mathbf{I}_{n-1}(\Rn)$ satisfies the conditions of Theorem \ref{thm:minimizing-weaklyreg}. Pick $x\in \spt T$ and $R$ such that $\Ball(x,R)\cap \spt \partial T=\emptyset$, then given $\nu>0$, there exists $C>0$ such that for $r\in (0,R)$, there holds
\[
    \Vert T\Vert (\Ball(x,R)\cap \Ball(\Sing(T),r)) \leq C r^{7-\nu}.
\]
\end{Theorem}

If $E_T= E_T^i$ then this is sufficient to prove Theorem \ref{thm:minimizing-weaklyreg}. However, there could a priori be singular points arbitrarily close to $\spt \partial T$. This was ruled out in the case where $\partial T$ represents a $C^{1,\alpha}$ submanifold by R. Hardt and L. Simon:
\begin{Theorem}[Boundary Regularity Theorem {\cite{HardtSimon1979}}]\label{thm:HardtSimon}
Under the assumptions of Theorem \ref{thm:minimizing-weaklyreg}, there exists a neighbourhood $V$ of $\spt \partial T$ such that $V  \cap \spt T$ is a $C^{1,\alpha}$ submanifold with boundary.
\end{Theorem}
Note that this theorem does not imply that $T$ has $C^1$-BV charts at a point of $\spt \partial T$. Indeed, there could be points of the boundary where $T$ has a jump in multiplicity when crossing $\spt \partial T$. However, combining the last two statements and a control of the Minkowski content of the boundary, we can obtain a proof of the main result of this paragraph.
\begin{proof}[Proof of Theorem \ref{thm:minimizing-weaklyreg}]
    Clearly $E_T \subset E_T^i \cup \spt \partial T$. By Theorem \ref{thm:HardtSimon}, $E_T^i$ is compact. As Minkowski content is additive on disjoint compact sets, it suffices to prove that $E_T^i$ and $\spt \partial T$ have finite Minkoswki content of dimension $m-1$ with respect to $\Vert T\Vert$. 

Covering $E_T^i$ by finitely many balls which do not intersect the boundary, we can apply Theorem \ref{thm:CheegerNaber} to obtain $C$ such that for for $r$ small enough, 
\[
 \Vert T\Vert (\Ball(\Sing(T),r)) \leq C r^6,
\]
this clearly implies that $\mathscr{M}^{m-1}_{\Vert T\Vert,*}(E^i_T) = 0$.

We will now control the Minkowski content of $\spt \partial T$. Given $x\in \spt \partial T$, by Theorem \ref{thm:HardtSimon} there exists $r>0$ such that $\spt T\cap \Ball(x,r)$ is an oriented $C^{1,\alpha}$ submanifold of dimension $n-1$ with boundary. Choosing $r$ smaller if necessary, we can suppose that $\spt \partial T$ separates $\spt T\cap \Ball(x,r)\backslash \partial T$ into at most two connected relatively open subsets of $\spt T$. We can therefore apply the constancy theorem for currents (see section 4.1.7 in \cite{FedererGMT}) to prove that the multiplicity of $T$ takes only one value in each connected component. Covering the compact set $\spt \partial T$ with finitely many such balls, we get an upper bound $\theta_{max}$ for the multiplicity of $\Vert T\Vert$ in a neighbourhood $V$ of $\spt \partial T$.
By standard differential geometry, there exists a constant $C>0$ such that for $r>0$ small enough, $\Ball(\spt \partial T,r)\subset V$ and 
\[
   \scrH^{n-1} (\spt T\cap \Ball(\spt \partial T,r)) \leq C \scrH^{n-2}(\spt \partial T)\ r,
\]
which implies
\[
\Vert T\Vert( \Ball(\spt \partial T,r)) \leq C \theta_{max} \scrH^{n-2}(\spt \partial T) \ r.
\]
Thus $\mathscr{M}^{m-1}_{\Vert T\Vert,*} (\spt \partial T) <+\infty$, which concludes the proof.
\end{proof}

\section{Semi-algebraic chains are weakly regular currents}\label{sec:semi-alg}
In this section we present another type of integral currents for which Theorem \ref{thm:stokes} holds. These currents can have singular sets of codimension $1$ and can represent for example real algebraic subvarieties of $\Rn$. We focus on semi-algebraic currents for simplicity, but we chose to give a proof which extends to any o-minimal structure. In \cite[Chapter 6]{ThesisAJ}, we gave a more detailed presentation. We refer the reader to \cite{VDD98book} for a comprehensive introduction to o-minimal geometry.

\begin{Definition}\label{def:semi-alg}
A set $E$ in $\Rn$ is \emph{semi-algebraic} if it can be written as a finite union of sets $E_k$ where for each $k$ there exists finitely many polynomials $P_k$ and $Q_k^1,\dots, Q_k^{q_k}$, such that
\begin{equation*}
    E_k= \{x\in \Rn, P_k(x)= 0, Q_k^1(x)>0, \dots, Q_k^{q_k}(x)>0\}.
\end{equation*}
\end{Definition}
In particular, the graphs and epigraphs of polynomials and rational fractions are semi-algebraic sets.

For $n=1,2,\dots$, denote by $\mathfrak{M}_n$ the collection of semi-algebraic subsets of $\Rn$; it satisfies the following properties
\begin{enumerate}[ref={*}]
\item For every $n$, $\mathfrak{M}_n$ is stable under union, intersection, and complement operations.
\item For every $n,n'$, $\mathfrak{M}_{n+n'}\supseteq \mathfrak{M}_n \times \mathfrak{M}_{n'} \defeq \left \{ A\times B, A\in \mathfrak{M}_n, B\in \mathfrak{M}_{n'}\right\}$.
\item For every $n$ and any $A\in \mathfrak{M}_n$, if $\pi:(x_1,\dots,x_n) \mapsto(x_1,\dots,x_{n-1})$, there holds $\pi(A)\in \mathfrak{M}_{n-1}$.
\end{enumerate}
The last point is highly non trivial, it is known as the Tarski-Seidenberg Theorem, see for instance \cite[1.4.6, 2.2.1]{BocCosRoy1998realalg}.

More generally, a \emph{structure on $\R$}, is a sequence for $n=1,2,\dots$ of collections $\mathfrak{M}_n\subset \calP(\Rn)$, containing the semi-algebraic sets and satisfying the three properties above. Such a structure is \emph{o-minimal} if $\mathfrak{M}_1$ consists of all the finite unions of intervals and points. Clearly the semi-algebraic sets form an o-minimal structure and we will stick to the semi-algebraic setting in the sequel. 

A map from a subset of $\Rn$ to $\Rm$ is semi-algebraic if its graph is a semi-algebraic set in $\R^{n+m}$. A first important property is that a semi-algebraic function from a subset of $\R$ to $\R$ is piecewise monotone and $C^1$. An important generalization of this property to higher dimensions is the so-called Cell Decomposition Theorem. Given a choice of coordinate axes, cells are defined iteratively as follows
\begin{enumerate}[label=(\roman*)]
  \item A $(0)$-cell in $\R$ is a point, a $(1)$-cell in $\R$ is an open interval.
  \item For $(i_1,\dots,i_{n-1})\in \{0, 1\}^{n-1}$, an $(i_1,\dots,i_{n-1}, 0)$-cell in $\R^{n}$ is the graph of a $C^1$ semi-algebraic function over an $(i_1,\dots,i_{n-1})$-cell in $\R^{n-1}$.
  \item an $(i_1,\dots,i_{n-1}, 1)$-cell in $\R^{m}$ is a set of the form
\[
    \{(x_1,\dots,x_n); (x_1,\dots, x_{n-1})\in A \quad \text{and} \quad (f(x_1,\dots, x_{n-1})<x_n<g(x_1,\dots, x_{n-1}))\},
\]
where $A$ is an $(i_1,\dots,i_{n-1})$-cell in $\R^{n-1}$ and $f$ and $g$ are $C^1$ semi-algebraic functions from $A$ to $\R$ with $f<g$.
\end{enumerate}
An \emph{$m$-cell} in $\Rn$ is an $(i_1,\dots,i_n)$-cell with $\sum_j i_j = m$. In particular, an $n$ cell in $\Rn$ is an open set. An $m$-cell in $\Rn$ is a $C^1$ graph over an $m$-cell in $\Rm$ up to a permutation of coordinate axes. However, the defining maps are not necessarily continuous up to the boundary, though if a cell is bounded, then the defining maps are bounded. Given an $m$-cell $A$, its frontier $\Fr(A)\defeq \cl A\backslash A$ is a finite union of cells of lower dimensions.

A \emph{cell decomposition} of $\R$ is a finite partition of $\R$ into points and open intervals. A \emph{cell decomposition} of $\Rn$ is a finite partition of $\Rn$ into cells which are all defined over a cell decomposition of $\R^{n-1}$. The importance of cell decompositions is underlined by the following result:
\begin{Theorem}[$C^1$ cell decomposition of semi-algebraic sets {\cite[Section~2.3]{BocCosRoy1998realalg}}]\label{thm:cell-decomp}
  Let $A_1,\dots, A_p$ be semi-algebraic sets in $\Rn$. Then there exists a cell decomposition $\calD$ of $\Rn$ such that each $A_j$ is a finite union of elements of $\calD$.
\end{Theorem}

The semi-algebraic setting is not a priori adapted to usual analytic tools. In particular taking limits of sequences of semi-algebraic objects is not interesting. However, there is a semi-algebraic way to take limits, via so called semi-algebraic families. If $B$ is a subset of $\Rm$ and for each $t\in B$, $A_t$ is a subset of $\Rn$, one says that $(A_t)_{t\in B}$ is a \emph{semi-algebraic family} if the set
\[
     A_B \defeq \bigcup_{t\in B} A_t\times \{t\} \subset \Rn\times \Rm
\]
is semi-algebraic. If $B$ is the interval $(0,1\,]$ and the $A_t$ are contained in the same compact set, one can define the limit of $A_t$ as $t$ goes to $0$ as $A_0 \defeq \cl (A_B)\cap \Rn\times \{0\}$. Then $A_0$ is a semi-algebraic set, whicht coincides with the Hausdorff limit of the sets $\cl A_t$.

Finally, while the derivative of a semi-algebraic function is a semi-algebraic function, this is not the case of its primitive --- take for instance the primitive of $x\mapsto 1/x$. For this reason, it is not easy to give measure estimates on semi-algebraic sets. However, the following result is sufficient for our purpose.
\begin{Theorem}[Uniform bound on the Hausdorff measure, {\cite{LoiPhi14}}]\label{thm:taleloi}
  If $(A_t)_{t\in B}$ is a bounded semi-algebraic family in $\Rn$, then for $m=0,\dots,n$, there exists a constant $C_m$ such that for $t\in B$, either $\scrH^m(A_t) <C_m$ or $A_t$ has dimension larger than $m$. The latter can happen only for a definable subset of indices in $B$.
\end{Theorem}

Let us now define a class of currents adapted to the representation of semi-algebraic objects. 
\begin{Definition}\label{def:semi-alg-chain}
A current in $\calD_m(\Rn)$ is a \emph{semi-algebraic cellular current} if it is of the form $\theta (\scrH^m\hel A) \wedge \vect T$, where $\theta$ is an integer, $A$ is a bounded semi-algebraic $C^1$ cell of dimension $1$ and $\vect T$ is a choice of orientation for $A$. A finite sum of semi-algebraic cellular currents of the same dimension $m$ is a \emph{semi-algebraic chain} of dimension $m$.
\end{Definition}

Using the Cell Decomposition Theorem \ref{thm:cell-decomp}, we can prove (see \cite[Proposition~6.3.1]{ThesisAJ}) that the cellular currents used in the definition of a semi-algebraic chain can be chosen according to the same cell decomposition. The semi-algebraic chain $T$ then admits a $C^1$ chart at all points of the defining cells. Using the last fact, cell decomposition and a variant of the Constancy Theorem \cite[4.1.7]{FedererGMT}, it can be proved that semi-algebraic chains are integral currents, whose boundary is a semi-algebraic chain (see \cite[Section~2]{Funk2016Homology} and \cite[Theorem~6.3.2]{ThesisAJ}). 

Notice that the singular set of a cellular current is contained in the support of its boundary, which is a semi-algebraic set of dimension $m-1$. Given a semi-algebraic chain $T$, we can chose an appropriate decomposition into cellular currents and infer that the singular set of $T$ is a semi-algebraic set of dimension $m-1$. We can actually say more:
\begin{Theorem}
   The singular set $E_T$ of a semi-algebraic chain $T$ is disposable. In particular $T$ is weakly regular.
\end{Theorem}
\begin{proof}
As $E_T$ is semi-algebraic, the map: $x\mapsto \dist(x,E_T)$ is a semi-algebraic map
 Thus, for $t\in [\,0,1\,]$, the set $A_t= \{x\in \spt T, \dist(x, E_T) = t\}$ is semi-algebraic, but more importantly the set
\[
     \{(t,x)\in [\,0,1\,] \times \Rn, x\in \spt T, \dist(x, \Sing T) = t\}
\]
forms a bounded definable family indexed by $t$. By Theorem \ref{thm:taleloi}, there exists $C$ such that for every $t$, $\scrH^{m-1}(E_t)$ is either less than $C$ or $E_t$ has dimension $m$. By classical slicing theory (or the coarea formula), the set of indices $t$ such that $E_t$ has dimension $m$ is Lebesgue null in $[\,0,\infty)$; it is also semi-algebraic and therefore must be discrete. Thus for small enough $t$, there holds $\scrH^{m-1}(E_t)<C$ and the current
\[
    T_t \defeq T\hel\{x, \dist(x,E_T)<t)\}
\]
is a semi-algebraic chain of dimension $m$, with
\[
\mass(\partial T_t) \leq \mass(\partial T)+  N \theta C,
\]
 where $N$ is the number of cellular currents involved in the construction of $T$ and $\theta$ is the maximal multiplicity of these currents. Furthermore, $\mass(T_t)$ tends to zero as $t$ goes to zero. Thus $E_T$ is disposable in $T$ and $T$ is weakly regular.
\end{proof}

\bibliographystyle{plain}
\bibliography{stokes}
\end{document}